\documentclass[12pt,a4paper]{article}
\usepackage[centertags]{amsmath}
\usepackage{amsfonts}
\usepackage{amssymb}
\usepackage{amsthm}
\usepackage{latexsym}
\usepackage{newlfont}
\usepackage{verbatim}
\usepackage[english]{babel}
\usepackage[arrow, matrix]{xy}
\usepackage[dvips]{graphicx}
\usepackage{rotating}
\newcommand{\Ae}[1]{#1^*}

\newcommand{\s}[1]{#1^*}
\newtheorem{thm}{Theorem}[section]
\newtheorem{cor}[thm]{Corollary}
\newtheorem{lem}[thm]{Lemma}
\newtheorem{prop}[thm]{Proposition}

\theoremstyle{definition}
\newtheorem{defn}{Definition}
\theoremstyle{remark}

\newcommand{\set}[1]{\left\{#1\right\}}
\newcommand{\nin}{\!\not\!\varepsilon\,}
\newcommand{\I}{\mathbb I}

\newcommand{\N}{\mathbb N}
\newcommand{\C}{\mathbb C}

\newcommand{\uf}{\mathtt{Uf}}
\def\bisim{{\underline{\leftrightarrow} }}
\def\proof{\trivlist \item[\hskip \labelsep{\bf Proof.\ }]}
\def\endproof{\null\hfill$\dashv$\endtrivlist}
\newenvironment{cproof}{\medskip\noindent{\sc Proof.\ }}
{\null\hfill\qedsymbol\endtrivlist}

\begin{document}

%{{{ abstract

\author{Balder ten Cate, David Gabelaia and Dmitry Sustretov}
\title{Modal Languages for Topology:
  Expressivity and Definability}

\maketitle

\begin{abstract}
  In this paper we study the expressive power and definability for (extended)
  modal languages interpreted on topological spaces. We provide topological analogues of the van
  Benthem characterization theorem and the Goldblatt-Thomason definability
  theorem in terms of the well established first-order topological language
  $\mathcal{L}_t$.
\end{abstract}

\pagebreak

{\small\tableofcontents}

\bigskip

\bigskip

\bigskip

{\large \bf Acknowledgements}

\medskip

Authors wish to thank Yde Venema for his invaluable help.

The work of the first author was supported by NWO grant 639.021.508,
and the work of the second author was supported by INTAS grant Nr
1-04-77-7080.
%}}}

\newpage
%{{{ Preliminaries

\section{Introduction}

Modal logic, as a language for talking about topological spaces, has
been studied for at least 60 years. Originally, the motivations for
this study were purely mathematical.  More recently, computer
science applications have lead to a revival of interest, giving rise
to new logics of space, many of which are (extensions of) modal
languages (e.g., \cite{Kremer-mints05, artemov97modal, kkwz06,
bbcs06}, to name a few).

The design of such logics is usually guided by considerations
involving expressive power and computational complexity. Within the
landscape of possible spatial languages, the basic modal language
interpreted on topological spaces can be considered a minimal
extreme. It has a low computational complexity, but also a limited
expressive power.

In this paper, we characterize the expressive power of the basic
modal language, as a language for talking about topological spaces,
by comparing it to the well established topological language
$\mathcal{L}_t$ \cite{flum-ziegler}. Among other things, we obtain
the following results:

\vspace{2mm}

\noindent {\bf Theorem \ref{t:vB charcterization}.} Let $\phi(x)$ be
any $\mathcal L_t$ formula with one free variable. Then $\phi(x)$ is
equivalent to (the standard translation of) a modal formula iff
$\phi(x)$ is invariant under topo-bisimulations.

\vspace{2mm}

\noindent {\bf Theorem \ref{t:GbTh basic Lt}.} Let $\mathsf K$ be a
class of topological spaces definable in $\mathcal L_t$. Then
$\mathsf K$ is definable in the basic modal language iff $\mathsf K$
is closed under topological sums, open subspaces and images of
interior maps, while the complement of $\mathsf K$ is closed under
Alexandroff extensions.

\vspace{2mm}

These can be seen as topological generalizations of the Van Benthem
theorem and the Goldblatt-Thomason theorem, respectively. We give
similar characterizations for some extensions of the modal language,
containing nominals, the global modality, the difference modality, and
the $\downarrow$-binder (for a summary of our main results, see
Section~\ref{s:discussion}).

Characterizations such as these help explain why certain languages (in
this case the basic modal language) are natural to consider. They can
also guide us in finding languages that provide the appropriate level
of expressivity for an application.

\subsubsection*{Outline of the paper}

The structure of the paper is as follows: Section
\ref{s:preliminaries} contains basic notions from topology,
topological model theory, and the topological semantics for modal
logic. Section~\ref{s:definability} is the core of the paper: in
Section~\ref{s:topo-bisimulations} we characterize the expressivity
of the basic modal language; Theorem~\ref{t:up_ae} of
Section~\ref{s:alexandroff ext} is the main technical result that is
used extensively in subsequent sections, while in Section~\ref{s:
GbTh main} we compare definability in the basic modal language with
first-order definability. Section~\ref{s:algebra} provides the
proper algebraic perspective on these results. In
Section~\ref{s:extensions}, we consider a number of extensions of
the basic modal language and characterize definability in these
richer languages. Finally, we conclude in
Section~\ref{s:discussion}.

\section{Preliminaries}\label{s:preliminaries}

In this section we recall some basic notions from topology,
topological model theory, and the topological semantics for
modal logic.

\subsection{Topological spaces}

\begin{defn}[Topological spaces]
A \emph{topological space} $(X,\tau)$ is a non-empty set $X$
together with a collection $\tau\subseteq \wp(X)$ of subsets that
contains $\emptyset$ and $X$ and is closed under finite
intersections and arbitrary unions. The members of $\tau$ are called
\emph{open sets} or simply \emph{opens}. We often use the same
letter to denote both the set and the topological space based on
this set: $X=(X,\tau)$.
\end{defn}

If $A\subseteq X$ is a subset of the space $X$, by $\I A$ (read:
`interior $A$') one denotes the greatest open contained in $A$
(i.e.\ the union of all the opens contained in $A$). Thus $\I$ is an
operator over the subsets of the space $X$. It is called the
\emph{interior operator}.

Complements of open sets are called \emph{closed}. The
\emph{closure} operator, which is a dual of the interior operator, is
defined as $\C A=-\I{-}A$ where `${-}$' stands for the set-theoretic
complementation. Observe that $\C A$ is the least closed set
containing $A$.

A standard example of a topological space is the real line
$\mathbb{R}$, where a set is considered to be open if it is a union
of open intervals $(a,b)$.

For technical reasons, at times it will be useful to consider
topological bases---collections of sets that generate a topology.

\begin{defn}[Topological bases] \label{def:topological-base} A
  \emph{topological base} $\sigma$ is a collection $\sigma\subseteq
  \wp(X)$ of subsets of a set $X$ such that closing $\sigma$ under
  arbitrary unions gives a topology on $X$ (i.e., such that
  $(X,\{\bigcup\sigma'\mid \sigma'\subseteq\sigma\})$ is a topological
  space). The latter requirement is in fact equivalent to the
  conjunction of the following conditions:
\begin{enumerate}
\item $\emptyset\in\sigma$
\item $\bigcup\sigma = X$
\item For all $A,B\in\sigma$ and $x\in A\cap B$, there is
 a $C\in\sigma$ such that $x\in C$ and $C\subseteq A\cap B$.
\end{enumerate}
For $(X,\sigma)$ a topological base, we denote by
$\widehat{X}=(X,\widehat{\sigma})$ the topological space it generates, i.e., the
topological space obtained by closing $\sigma$ under arbitrary
unions. Furthermore, we say that $\sigma$ \emph{is a base for}
$\widehat{\sigma}$.
\end{defn}

For example, a base for the standard topology on the reals is the
set of open intervals $\set{(a,b)\mid a\leq b}$.

\subsection{The basic modal language}\label{s:modal}

We recall syntax and the topological semantics for the basic modal
language.

\begin{defn}[The basic modal language]
  The \emph{basic modal language} $\mathcal{ML}$ consists of a set of
  propositional letters $\textsc{Prop}=\{p_1, p_2, \dots\}$, the
  boolean connectives $\wedge, \neg$, the constant truth $\top$ and a
  modal box $\Box$.  Modal formulas are built according to the
  following
  recursive scheme:\\
$$\phi\ \ ::=\  \top ~\mid~ p_i ~\mid~ \phi\wedge\phi ~\mid~ \neg\phi ~\mid~ \Box\phi$$
\end{defn}
We use $\Diamond \phi$ as an abbreviation for $\neg \Box \neg \phi$.
Unless specifically indicated otherwise, we will always assume that
the set of propositional letters is countably infinite.

Nowadays, the best-known semantics for $\mathcal{ML}$ is the Kripke
semantics. In this paper, however, we study the topological
semantics, according to which modal formulas denote regions in a
topological space. The regions denoted by the propositional letters
are specified in advance by means of a \emph{valuation}, and
$\land$, $\neg$ and $\Box$ are interpreted as intersection,
complementation and the interior operator.  Formally:

\begin{defn}[Topological models]
  A {\em topological model} $\mathfrak{M}$ is a tuple $(X,\nu)$ where
  $X=(X,\tau)$ is a topological space and the valuation
  $\nu:\textsc{Prop}\to\wp(X)$ sends propositional letters to subsets
  of $X$.
\end{defn}

\begin{defn}[Topological semantics of the basic modal language]
  \label{d:topo-semantics}
  Truth of a formula $\phi$ at a point $w$ in a topological
  model $\mathfrak{M}$ (denoted by $\mathfrak{M},w \models \phi$) is
  defined inductively:
\[\begin{array}{lll}
\mathfrak{M},w \models \top & & \mathrm{always}\\
\mathfrak{M},w \models p & \mathrm{\mathrm{iff}} & x \in \nu(p)\\
\mathfrak{M},w \models \phi \land \psi & \mathrm{iff} &
\mathfrak{M},w
\models \phi \textrm{ and } \mathfrak{M},w \models \psi\\
\mathfrak{M},w \models \neg \phi & \mathrm{iff} & \mathfrak{M},w
\nvDash
\phi\\
\mathfrak{M},w \models \Box \phi & \mathrm{iff} & \exists O \in \tau
\textrm{ such that } w\in O \textrm{ and } \forall v\in
O.(\mathfrak{M},v \models \phi)
\end{array}\]
  If $\mathfrak{M},w \models \phi$ for all $w \in A$ for some $A\subset
  X$, we write $A \models \phi$. Further,
  $\mathfrak{M} \models \phi$ ($\phi$ is \emph{valid} in $\mathfrak{M}$)
  means that $\mathfrak{M},w \models \phi$ for all $w \in X$. We write
  $X \models \phi$ ($\phi$ is valid in $X$) when $(X,\nu) \models \phi$
  for any valuation $\nu$. If $\mathsf K$ is a class of topological
  spaces we write $\mathsf{K} \models \phi$ when $X \models \phi$ for
  each $X \in \mathsf{K}$.
\end{defn}
Each modal formula $\phi$ defines a set of points in a topological
model (namely the set of points at which it is true).  With a slight
overloading of notation, we will sometimes denote this set by
$\nu(\phi)$. It is not hard to see  that $\nu(\Box\phi)=\I\nu(\phi)$.

We extend the notions of truth and validity to the sets of
modal formulas in the usual way (e.g.,
$X\models \Gamma$ means that $X\models\phi$ for each
$\phi\in\Gamma$).

\begin{defn}[Modal definability]
  A set of modal formulas $\Gamma$ \emph{defines} a class $\mathsf{K}$
  of spaces if, for any space $X$, $$X \in \mathsf{K}\mathrm{\ \ \ iff\ \ \ }X
  \models \Gamma$$  A class of topological spaces is said to be
  \emph{modally definable} if there exists a set of modal formulas that defines it.
A topological property is said to be modally definable
  (or, defined by a set of formulas $\Gamma$) if the class of all
  spaces that have the property is modally definable (is defined by
  $\Gamma$).
\end{defn}

  Given a class $\mathsf K$ of spaces, the set of modal formulas
  $\set{\phi\in \mathcal{ML}\mid \mathsf K\models\phi}$ (``\emph{the
    modal logic of $\mathsf{K}$}'') is denoted by
  $\mathrm{Log}(\mathsf K)$.  Conversely, given a set of modal formulas
  $\Gamma$, the class of spaces $\set{X\mid X\models \Gamma}$
  is denoted by $\mathrm{Sp}(\Gamma)$.  Thus, in this notation, a class
  $\mathsf K$ is modally definable iff $\mathrm{Sp}
  (\mathrm{Log}(\mathsf K))=\mathsf K$.

  The following example illustrates the concept of modal definability.

\begin{defn}[Hereditary Irresolvability]\label{def:irresolvability}
  A subset $A\subseteq X$ of a space $X$ is said to be \emph{dense} in
  $X$ if $\C A=X$ (or, equivalently, if $A$ intersects each non-empty
  open in $X$).
  A topological space $X$ is called \emph{irresolvable} if it cannot
  be decomposed into two disjoint dense subsets.  It is
  \emph{hereditarily irresolvable} (HI) if all its subspaces\footnote{Recall that a \emph{subspace} of a space $X$ is a
non-empty subset $A\subseteq X$ endowed with the relative topology
$\tau_A=\set{O\cap A\mid O\in\tau}$.} are
  irresolvable.
\end{defn}

\begin{thm}\label{t:HI definable}
  The modal formula $\Box(\Box(p\to\Box p)\to p)\to \Box p$ (Grz)
  defines the class of hereditarily irresolvable spaces.
\end{thm}

\proof Follows from results in \cite{Es81} and \cite{Gu03}.  For
purposes of illustration, we will give a direct proof, inspired by
\cite{Gu03}.

We are to show that $X$ is HI iff $X\models(Grz)$.

First note that $X\models \text{(Grz)}$ iff $X\models \Diamond \neg
p\to \Diamond(\neg p\wedge\Box(p\to\Box p))$ iff $X\models \Diamond
q\to \Diamond(q\wedge\neg\Diamond(\neg q\wedge \Diamond q))$ iff
$\forall A\subseteq X.[\C A\subseteq \C (A-\C(\C A-A))]$.

Suppose $X$ is not HI. Then there exists a non-empty subset  $A\subseteq X$
and two disjoint sets $B,B'\subset A$ such that $A\subseteq\C B\cap \C
B'$. We show that $\C B\not\subseteq \C (B-\C(\C B-B))$ so $X$ does not make $(Grz)$ valid.
Indeed, since $A\subseteq \C B$ it is clear that $B'\subseteq \C B - B$, hence $B\subseteq A\subseteq \C B'\subseteq \C(\C
B-B)$ and $\C B\not\subseteq \C (B-\C(\C B-B))=\emptyset$.

Suppose $X\not\models (Grz)$. Then there exists a non-empty subset
$A\subseteq X$ such that $\C A\not\subseteq \C (A-\C(\C A-A))$.
Denote $Y=\C A$. We will show that $Y$ is not HI thus proving that
$X$ is not HI (it is easily seen that a closed subspace of an HI
space must itself be HI). Since $Y$ is a closed subspace of $X$ the
operator $\C_Y$ coincides with $\C$ on subsets of $Y$. Thus
$Y\not\subseteq \C_Y (A-\C_Y(Y-A))=\C_Y\I_Y A$. It follows that $A$
is dense in $Y$ while $\I_Y A$ is not dense in $Y$. Then there
exists a subset $U\subseteq Y$ that is open in the relative topology
of $Y$ such that $\emptyset=\I_Y A\cap U=\I_Y(U\cap
A)=\I_Y((U\cap(Y-U))\cup (U\cap A))=$ $\I_Y(U\cap(A\cup(Y-U)))=
U\cap\I_Y(A\cup(Y-U))=U-\C_Y(U-A)$. This implies that $U\subseteq
\C_Y(U-A)$. But at the same time $U\subseteq \C_Y(U\cap A)$ since
$U$ is open in $Y$ and $A$ is dense in $Y$. As $U=(U-A)\cup(U\cap
A)$ it follows that $U$ is decomposed into two disjoint dense in $U$
subsets $U-A$ and $U\cap A$, so $U$ is resolvable. Thus $Y$ is not
HI and hence $X$ is not HI either.
\endproof

% \textbf{\sc Question:} Is this property elementary? My guess is that
% it is expressible in $\mathcal L^2$, but not in $\mathcal L_t$.
% This is unfortunate because Sahlqvist fails...\\

One of the central questions in this paper is which properties of
topological spaces are definable in the basic modal language and its
various extensions.

\subsection{The topological correspondence language
  $\mathcal{L}_t$}\label{s:elementarity}

In the relational semantics, the van Benthem theorem and the
Goldblatt-Thomason theorem characterize the expressive power of the
basic modal language by comparing it to the `golden standard' of
first-order logic. In the topological setting, it is less clear what
the golden standard should be.  Let us imagine for a moment a perfect
candidate for a `first-order correspondence language for topological
semantics of modal logic'.  Such a language should have the usual kit
of nice properties of first-order languages like Compactness and the
L\"{o}wenheim-Skolem theorem; it should be able to express topological
properties in a natural way; moreover, it should be close enough to
the usual mathematical language used for speaking about topologies so
that we could determine easily whether a given topological property is
expressible in it or not; and it should be suitable for translating modal
formulas into it nicely.

The language $\mathcal L_t$ which we describe in this section
satisfies all these requirements. Moreover, its model theory has been
quite well studied and the corresponding machinery will serve us well
in the following sections. With the exception of
Theorems~\ref{thm:L2-is-bad} and \ref{thm:Lt-saturation}, all results on $\mathcal{L}_t$
discussed in this section, and much more, can be found in the
classical monograph on topological model theory by Flum and Ziegler
\cite{flum-ziegler}.

Before defining $\mathcal L_t$, we will first introduce the two-sorted
first order language $\mathcal L^2$. In its usual definition, this
language can contain predicate symbols of arbitrary arity. Here,
however, since the models we intend to describe are the
\emph{topological models} introduced in the previous section, we will
restrict attention to a specific signature, containing a unary
predicate for each propositional letter $p\in\textsc{Prop}$.

\begin{defn}[The quantified topological language $\mathcal L^2$]
  $\mathcal L^2$ is a two-sorted first-order language: it has terms that are
  intended to range over elements, and terms that are intended to range
  over open sets. Formally, the alphabet is constituted by a countably infinite set of ``point
  variables'' $x, y, z, \ldots$ a countably infinite set of ``open
  variables'' $U, V, W, \ldots$, unary predicate symbols $P_p$ corresponding to
   propositional letters $p\in\textsc{Prop}$ and a binary predicate
   symbol $\varepsilon$ that relates point variables with open variables. The formulas of $\mathcal L^2$ are given by
  the following recursive definition:
\[ \phi ~::=~ \top ~\mid~ x=y ~\mid~ U=V ~\mid~ P_p(x) ~\mid~ x\varepsilon U ~\mid~ \neg\phi ~\mid~
   \phi\land\phi ~\mid~ \exists x.\phi ~\mid~ \exists U.\phi \]
   where $x,y$ are point variables and $U,V$ are open
   variables. The usual shorthand notations (e.g., $\forall$ for
   $\neg\exists\neg$) apply.
\end{defn}
Due to the chosen signature, formulas of $\mathcal{L}^2$ can be
naturally interpreted in topological models (relative to assignments
that send point variables to elements of the domain and open
variables to open sets). However, as we show in
Appendix~\ref{app:L2}, under this semantics, $\mathcal{L}^2$ is
rather ill-behaved: it lacks the usual model theoretic features such
as Compactness, the L\"owenheim-Skolem theorem and the \L{}o\'s
theorem.  For this reason, we will first consider a more general
semantics in terms of \emph{basoid
  models}.

\begin{defn}[Basoid models]
  A \emph{basoid model} is a tuple $(X,\sigma, \nu)$ where $X$ is a non-empty set,
  $\sigma\subseteq\wp(X)$ is a topological base, and the valuation $\nu:\textsc{Prop}\to\wp(X)$
  sends propositional letters to subsets of $X$.
\end{defn}

Interpret $\mathcal L^2$ on a basoid model as follows: point
variables range over $X$, open variables range over $\sigma$, the
valuation $\nu$ determines the meaning of the unary predicates
$P_p$, while $\varepsilon$ is interpreted as the set-theoretic
membership relation.

Under this interpretation, $\mathcal L^2$ displays all the usual
features of a first-order language, including Compactness, the
L\"{o}wenheim-Skolem property and the \L{}o\'s theorem  \cite{flum-ziegler}.%
\footnotemark  As we mentioned already, these properties are lost if
we further restrict attention to topological models.

\footnotetext{Essentially, this is due to the fact that, within the
class
  of all two-sorted first-order structures, the basoid models can be
  defined up to isomorphism by conjunction of the following sentences of $\mathcal{L}^2$ (cf.~Definition~\ref{def:topological-base}, see also \cite[p.
  14]{garavaglia}):
$$\begin{array}{rcl}
    \textrm{Ext} & \equiv & \forall U,V.(U=V\leftrightarrow
        \forall x.(x\varepsilon U\leftrightarrow x\varepsilon V)) \\
    \textrm{Union} & \equiv &\forall x.\exists U.(x\varepsilon U)\\
    \textrm{Empty} & \equiv &\exists U.\forall x.(\neg x\varepsilon U)\\
    \textrm{Bas} & \equiv &\forall x.\forall U,V.(x\varepsilon U\wedge x\varepsilon V\to
        \exists W.(x\varepsilon W\wedge\forall z.(z\varepsilon W\to
        z\varepsilon U\wedge z\varepsilon V)))
\end{array}$$
}

\begin{thm}\label{thm:L2-is-bad}
  $\mathcal{L}^2$ interpreted on topological models
  lacks Compactness, L\"owenheim-Skolem and Interpolation, and
  is $\Pi^1_1$-hard for validity.
\end{thm}

The proof can be found in Appendix~\ref{app:L2}.

Thus, in order to work with topological models \emph{and} keep the nice
first-order properties we need to somehow `tame' $\mathcal L^2$. This
is where $\mathcal L_t$ enters the picture, a well behaved fragment of $\mathcal
L^2$.  Let us call an $\mathcal L^2$ formula $\alpha$ {\em positive
  (negative) in an open variable $U$} if all free occurrences of $U$
are under an even (odd) number of negation signs.

\begin{defn}[The language $\mathcal{L}_t$]
$\mathcal  L_t$ contains all atomic $\mathcal L^2$-formulas and is
closed under conjunction, negation, quantification over the point
variables and the following restricted form of quantification over
open variables:
\begin{itemize}
\item[-] if $\alpha$ is positive in the open variable $U$, and $x$ is
a point variable, then $\forall U.(x \varepsilon U \to \alpha)$ is a
formula of $\mathcal L_t$,
\item[-] if $\alpha$ is negative in the open variable $U$, and $x$ is
a point variable, then $\exists U.(x \varepsilon U \land \alpha)$ is a
formula of $\mathcal L_t$.
\end{itemize}
 (recall that $\phi \to \psi$ is simply an abbreviation for $\neg
(\phi\wedge\neg\psi)$).
\end{defn}

The reason $\mathcal L_t$ is particularly well-suited for describing
topological models lies in the following observation: $\mathcal
L_t$-formulas cannot distinguish between a basoid model and the
topological model it generates. More precisely, for any basoid model
$\mathfrak{M}=(X,\sigma,\nu)$, let $\widehat{\mathfrak{M}} =
(X,\widehat{\sigma},\nu)$, where $\widehat{\sigma}$ is the topology
generated by the topological base $\sigma$.

\begin{thm}[$\mathcal{L}_t$ is the base-invariant fragment of $\mathcal{L}^2$]
  \label{thm:Lt-invariance}
  For any $\mathcal{L}_t$-formula $\alpha(x_1, \ldots, x_n,U_1,
  \ldots, U_m)$, basoid model $\mathfrak{M}=(X,\sigma,\nu)$, and for
  all $d_1, \ldots, d_n\in X$ and $O_1, \ldots, O_m\in\sigma$,
  \[
  \mathfrak{M}\models\alpha\,[d_1, \ldots, d_n,O_1, \ldots, O_m] \text{ ~iff~ }
  \widehat{\mathfrak{M}}\models\alpha\,[d_1, \ldots, d_n, O_1, \ldots,
  O_m]~.
  \]
  Moreover, every $\mathcal{L}^2$-formula $\phi(x_1, \ldots,
  x_n, U_1, \ldots, U_m)$ satisfying this invariance property is
  equivalent on topological models to
  an $\mathcal{L}_t$-formula with the same free variables.
\end{thm}

%As an example, consider the formula
%$\textrm{Bas}$. It is clearly invariant for topologies and indeed it
%has the following $\mathcal L_t$-equivalent:
%
%\[ \forall x.\forall U.(x\varepsilon U\to\forall V.(x\varepsilon V\to
%\exists W.(x\varepsilon W\wedge\forall z.(z\varepsilon W\to
%z\varepsilon U\wedge z\varepsilon V))))\]
%

It follows that $\mathcal L_t$ satisfies appropriate analogues of
Compactness, the L\"{o}wenheim-Skolem property, and the \L{}o\'s
theorem \emph{relative to the class of topological models}.
Let us start with the L\"{o}wenheim-Skolem property.  Call a
topological model $\mathfrak{M}=(X,\tau,\nu)$ \emph{countable} if $X$
is countable and $\tau$ has a countable base.
%Then we have the
%following analogue of the L\"owenheim-Skolem theorem for
%$\mathcal{L}_t$:

\begin{thm}[L\"owenheim-Skolem for $\mathcal{L}_t$] \label{thm:LoSk-Lt}
  Let $\Gamma$ be any set of $\mathcal{L}_t$-formulas (in a countable
  signature). If $\Gamma$ has an infinite topological model, then it has a
  countable topological model.
\end{thm}

%\proof Let $g$ be an assignment s.t. $\mathfrak{M}\models\Gamma[g]$.
%By the L\"owenheim-Skolem property of $\mathcal{L}^2$ relative to
%the class of basoid models, there is a countable basoid submodel
%$\mathfrak{N}$ of $\mathfrak{M}$ such that
%$\mathfrak{N}\models\Gamma[g]$. Substituting the topological model
%$\widehat{\mathfrak{N}}$ for $\mathfrak N$ affects neither the
%cardinality of the set of points, nor the base cardinality, nor the
%truth of the formulas in $\Gamma$.\endproof

Next, we will discuss an analogue of \L{}o\'s's theorem for
$\mathcal{L}_t$. First we need to define ultraproducts of topological
models.

\begin{defn}[Ultraproducts of basoid models]
  Let $(\mathfrak{M}_i)_{i\in I}$ be an indexed family of basoid models,
  where  $\mathfrak{M}_i=(X_i,\sigma_i,\nu_i)$, and let $\mathfrak D$ be an
  ultrafilter over the index set $I$. Define an equivalence relation
  $\sim_{\mathfrak D}$ on $\prod_{i \in I} X_i$ as follows:
  $$x \sim_{\mathfrak D} y\ \ \ \  \textrm{iff}\ \ \ \ \{i \mid x_i = y_i\} \in \mathfrak D$$
  We define the
  \emph{ultraproduct} $\prod_{\mathfrak D} \mathfrak{M}_i$ to be $(X,\sigma,\nu)$,
  where $X=(\prod_{i\in I} X_i)/_{\sim_{\mathfrak D}}$, $\sigma = \{(\prod_{i\in
    I}O_i)/_{\sim_{\mathfrak D}} \mid \text{each } O_i\in \sigma_i\}$, and $\nu(p) =
  (\prod_{i\in I}\nu_i(p))/_{\sim_{\mathfrak D}}$.

  If $\mathfrak{M}_i = \mathfrak{M_j}$ for all $i,j\in I$, then
  $\prod_{\mathfrak D} \mathfrak{M}_i$ is called an \emph{ultrapower}.
\end{defn}
It is not hard to see that, under this definition, every ultraproduct
of basoid models is again a basoid model. The same does \emph{not}
hold for topological models. Hence, rather than the basoid ultrapower
$\prod_{\mathfrak D} \mathfrak{M}_i$, we will use the topological
model it generates, i.e., $\widehat{\prod_{\mathfrak D}
  \mathfrak{M}_i}$.  We will call the latter the \emph{topological
  ultraproduct} (or, \emph{topological ultrapower}, if all factor
models coincide).  Note that, by Theorem~\ref{thm:Lt-invariance},
the topological ultraproduct $\widehat{\prod_{\mathfrak D}
  \mathfrak{M}_i}$ cannot be distinguished from the basoid
ultraproduct $\prod_{\mathfrak D} \mathfrak{M}_i$ in $\mathcal{L}_t$.

\begin{thm}[\L{}o\'s theorem for $\mathcal{L}_t$] \label{t:Los-Lt}
  Let $\alpha$ be any $\mathcal{L}_t$-sentence, $(\mathfrak M_i)_{i\in I}$ an
  indexed set of topological models, and $\mathfrak D$ an ultrafilter over $I$.
  Then $$\widehat{\prod_{\mathfrak D} \mathfrak M_i}\models\alpha\ \ \ \textrm{iff}\ \ \  \{i\in I\mid
  \mathfrak M_i\models\alpha\}\in \mathfrak D$$
  In particular, if $\mathfrak{N}$ is a
  topological ultrapower of $\mathfrak{M}$, then for all
  $\mathcal{L}_t$-formulas $\phi$ and assignments $g$,
  $\mathfrak{M}\models\phi~[g]$ iff $\mathfrak{N}\models\phi~[f\cdot g]$,
  where $f:\mathfrak{M}\to\mathfrak{N}$ is the natural diagonal embedding.
\end{thm}
A typical use of ultraproducts is for proving compactness.

\begin{thm}[Compactness for $\mathcal{L}_t$]\label{thm:compactness}
  Let $\Gamma$ be any set of $\mathcal{L}_t$-formulas.  If every
  finite subset of $\Gamma$ is satisfiable in a topological model,
  then $\Gamma$ itself is satisfiable in a topological ultraproduct
  of these models.
\end{thm}
%
%\proof Follows from compactness via ultraproducts for
%$\mathcal{L}^2$ (with respect to basoid models) using
%Theorem~\ref{thm:Lt-invariance}.  \endproof
%
Another common use of ultraproducts is for obtaining \emph{saturated
  models}.  One can generalize this construction to topological
models, provided that the notion of saturation is defined carefully
enough.  The following definition of saturatedness is probably not
the most general, but will suffice for the purposes of this paper.

\begin{defn}[$\mathcal{L}_t$-saturatedness]\label{def:Lt-saturation}
  By an \emph{$\mathcal{L}_t$-type} we will mean a set of $\mathcal{L}_t$-formulas
  $\Gamma(x)$ having exactly one free point variable $x$ and no free
  open variables. An open set $O$ in a topological model is called
  \emph{point-saturated} if, whenever all finite subtypes of an
  $\mathcal{L}_t$-type $\Gamma(x)$ are realized somewhere in $O$,
  then $\Gamma(x)$ itself is realized somewhere in $O$.
  A topological model $\mathfrak{M}=(X,\tau,\nu)$ is said to
  be \emph{$\mathcal{L}_t$-saturated} if the following conditions
  hold:
  \begin{enumerate}
  \item The entire space $X$ is point-saturated.

  \item The collection of all point-saturated open sets forms a base
    for the topology. Equivalently, for each point $d$ with open
    neighborhood $O$, there is a point-saturated open subneighborhood
    $O'\subseteq O$ of $d$.

  \item Every point $d$ has an open neighborhood $O_d$ such that, for
    all $\mathcal{L}_t$-formulas $\phi(x)$, if $\phi(x)$ holds
    throughout \emph{some} open neighborhood of $d$ then $\phi(x)$
    holds throughout $O_d$.
  \end{enumerate}
\end{defn}

\begin{thm} \label{thm:Lt-saturation}
  Every topological model $\mathfrak{M}$ has an
  $\mathcal{L}_t$-saturated topological ultrapower. This
  holds regardless of the cardinality of the language.
\end{thm}

\begin{proof}
Let $\mathfrak{M}$ be any topological model. It follows from
classical model theoretic results that $\mathfrak{M}$ has a basoid
ultrapower $\prod_{\mathfrak D} \mathfrak{M} = (X,\sigma,\nu)$ that
is countably saturated (in the classical sense, for the
language $\mathcal{L}^2$) \cite[Theorem 6.1.4 and
6.1.8]{chang-keisler}.  We claim that $\widehat{\prod_{\mathfrak D}
\mathfrak{M}}$ is $\mathcal{L}_t$-saturated.

In what follows, with \emph{basic open sets} we will mean
open sets from the basoid model $\prod_{\mathfrak D} \mathfrak{M}$.

\begin{enumerate}
\item Suppose every finite subset of an $\mathcal{L}_t$-type $\Gamma(x)$ is satisfied by some
  point in $\widehat{\prod_{\mathfrak D} \mathfrak{M}}$. In other
  words, for every finite $\Gamma'(x)\subseteq \Gamma$,
  $\widehat{\prod_{\mathfrak D} \mathfrak{M}}\models\exists
  x.\bigwedge\Gamma'(x)$. Note that the latter formula belongs to
  $\mathcal{L}_t$. It follows by base invariance (Theorem~\ref{thm:Lt-invariance}) that every finite
  subset of $\Gamma(x)$ is satisfied by some point in
  $\prod_{\mathfrak D} \mathfrak{M}$. Hence, by the countable
  saturatedness of this basoid model,
  there is a point $d$ satisfying all formulas of $\Gamma(x)$.
  Applying the base invariance again, we conclude that
  $d$ still satisfies all formulas of $\Gamma(x)$ in
  $\widehat{\prod_{\mathfrak D} \mathfrak{M}}$.

\item Let $d$ be any point and $O$ any open neighborhood of $d$. By
  definition, $O$ is a union of basic open sets from $\prod_{\mathfrak
    D} \mathfrak{M}$. It follows that $d$ must have a basic open
  subneighborhood $O'$. Of course, $O'$ is still an open neighborhood
  of $d$ in $\widehat{\prod_{\mathfrak D} \mathfrak{M}}$.  By the same
  argument as before, we know that $O'$ is point-saturated---just consider the type $\Sigma(x)=\set{x\varepsilon O'}\cup\Gamma(x)$.

\item Let $d$ be any point and let $\Sigma$ be the collection of all
  $\mathcal{L}_t$-formulas $\phi(x)$ that hold throughout some open
  neighborhood of $d$. Recall that each open neighborhood of $d$
  contains a basic open subneighborhood of $d$. It follows that each
  $\phi(x)\in\Sigma$ holds throughout some basic open neighborhood
  of $d$.

  Next, we will proceed using the language $\mathcal{L}^2$, and the
  fact that $\prod_{\mathfrak D} \mathfrak{M}$ is countably saturated
  as a model for this language.  Consider the following set of
  $\mathcal{L}^2$-formulas (where $\mathbf{d}$ is used as a parameter
  referring to $d$, and $U$ is a free open variable):
  $$\Gamma(U) = \{ \mathbf{d}\varepsilon U\} \cup
  \{ \forall y.(y\varepsilon U\to \phi(y)) \mid \phi(y)\in\Sigma\}$$
  Every finite subset of $\Gamma(U)$ holds throughout some
  basic open neighborhood of $d$ (in $\prod_{\mathfrak D}
  \mathfrak{M}$). This follows from the definition of $\Sigma$, the
  base invariance of $\mathcal{L}_t$, and the fact that every open
  neighborhood of $d$ contains a basic open neighborhood.

  It follows by the countable saturatedness of $\prod_{\mathfrak D}
  \mathfrak{M}$ with respect to $\mathcal{L}^2$ that there is a basic
  open set $O_d$ satisfying all formulas in $\Gamma(U)$. In particular,
  $O_d$ is an open neighborhood of $d$ and (applying base
  invariance once more) all formulas in $\Sigma$ hold
  throughout $O_d$ in $\widehat{\prod_{\mathfrak D} \mathfrak{M}}$.
\end{enumerate}
\vspace{-9mm}
\end{proof}

We can conclude that $\mathcal{L}_t$ is model theoretically quite well
behaved. Computationally, $\mathcal{L}_t$ is unfortunately less well
behaved.

\begin{thm}The $\mathcal{L}_t$-theory of all topological spaces is
  undecidable, even in the absence of unary predicates. The same holds
  for $T_0$-spaces, for $T_1$-spaces, and for $T_2$-spaces. The
  $\mathcal{L}_t$-theory of topological models based on $T_3$-spaces,
  on the other hand, is decidable.
\end{thm}

The next natural question is which topologically interesting
properties we can express in this language.
Table~\ref{tab:Lt-examples} lists some examples of properties that
can be expressed in $\mathcal{L}_t$ (where $x\nin U$ is used as
shorthand for $\neg(x\varepsilon U)$, $\forall U_x.\alpha$ as
shorthand for $\forall U.(x\varepsilon U\to\alpha)$, and $\exists
U_x.\alpha$ as shorthand for $\exists U.(x\varepsilon
U\wedge\alpha)$). Recall that the \emph{separation axioms} $T_0-T_5$
are properties of spaces that allow separating distinct points
and/or disjoint closed sets \cite{engelking}.

\begin{table}[h]
\hrulefill
\caption{Some examples of properties that can be expressed in
$\mathcal{L}_t$} \label{tab:Lt-examples}
\[\begin{array}{ll}
  T_0 & \forall xy.(x\neq y\to \exists U_x.(y\nin U)\vee \exists V_y.(x\nin V)) \\[0.2em]
  T_1 & \forall xy.(x\neq y\to \exists U_x.(y\nin U)) \\[0.2em]
  T_2 & \forall xy.(x\neq y\to \exists U_x.\exists V_y.\forall z.(z\nin U\vee z\nin V)) \\[0.2em]
  \text{Regular} &
  \forall x.\forall U_x.\exists V_x.\forall y.(y\varepsilon U\lor
    \exists V'_y.\forall z.(z\varepsilon V'\to z\nin V)) \\[0.2em]
  T_3 & T_2\land\text{Regular} \\[0.2em]
  \text{Discrete} & \forall x.\exists U_x.\forall y.(y\varepsilon U\to y=x) \\[0.2em]
  \text{Alexandroff} &
  \forall x.\exists U_x.\forall V_x.\forall y.(y\varepsilon V\to y\varepsilon U)
\end{array}\]
\hrulefill
\end{table}

A typical example of a property not expressible in $\mathcal{L}_t$
is \emph{connectedness}.

\begin{defn}[Connectedness] \label{def:connectedness}
  A topological space $(X,\tau)$ is said to
  be \emph{connected} if $\emptyset$ and $X$ are the only sets that
  are both open and closed.
\end{defn}

\begin{thm}[{\cite[page 8]{flum-ziegler}}] \label{thm:connectedness-Lt}
  Connectedness is not expressible in $\mathcal{L}_t$.
\end{thm}

Note that connectedness \emph{is} expressible in $\mathcal{L}^2$,
namely by the sentence $\forall U,U'.(\forall x.(x\varepsilon
U\leftrightarrow x\nin U')\to (\forall x.(x\varepsilon U) \lor
\forall x.(x\nin U)))$.

We have the following translation from the basic modal language to
$\mathcal{L}_t$.\footnote{In fact, a slight variation of
  this translation shows that modal formulas can be mapped to
  $\mathcal{L}_t$-formulas containing at most two point variables and
  one open variable.}
\begin{defn}\label{d:ST}[Standard translation] The \emph{standard translation} $ST$ from the basic modal
language $\mathcal{ML}$
  into $\mathcal L_t$ is defined inductively:
\[\begin{array}{lll}
ST_x(\top) & = & \top\\
ST_x(p) & = & P_p(x)\\
ST_x(\neg \phi) & = & \neg ST_x(\phi)\\
ST_x(\phi \land \psi) & = & ST_x(\phi) \land ST_x(\psi)\\
ST_x(\Box\phi) & = & \exists U. (x \varepsilon U \land \forall y.(y
\varepsilon U \to ST_y(\phi)))
\end{array}\]
where $x, y$ are distinct point variables and $U$ is an open variable.
\end{defn}

\begin{thm}\label{t:stand-transl}
For ${\mathfrak M}$ a topological model and  $\varphi \in
\mathcal{ML}$ a modal formula,
${\mathfrak M}, a\models\varphi$ iff ${\mathfrak
  M}\models ST_x(\varphi)[a]$
\end{thm}
\proof By induction on the complexity of $\varphi$.\endproof

%In what follows by an \emph{elementary class of spaces} we will mean
%a class of topological spaces definable by a set of $\mathcal
%L_t$-sentences. A class of spaces will be called  \emph{basic
%elementary}, if it is definable by a single $\mathcal L_t$-sentence.

In other words, modal formulas can be seen as $\mathcal{L}_t$-formulas
in one free variable, and sets of modal formulas can be seen as
\emph{$\mathcal{L}_t$-types} in the sense of
Definition~\ref{def:Lt-saturation}. This shows that all the above
results on $\mathcal{L}_t$ also apply to modal formulas. For example,

\begin{thm}[L\"{o}wenheim-Skolem theorem for $\mathcal{ML}$]
Let $\Sigma\subseteq\mathcal{ML}$ be a set of modal formulas (in a
countable signature).
If $\Sigma$ is satisfied in a topological model, then it
is satisfied in a countable topological model.
\end{thm}
%\proof Suppose
%$\mathfrak{M},w\models\varphi$ for all $\phi\in\Sigma$. Then, by
%Theorem~\ref{t:stand-transl}, $\mathfrak{M}\models ST_x(\varphi)~[w]$
%for all $\varphi\in\Sigma$. By Theorem~\ref{thm:LoSk-Lt},
%there is a countable topological model $\mathfrak{M}'$ with
%a point $w'$ such that $\mathfrak{M}'\models ST_x(\varphi)~[w']$
%for all $\varphi\in\Sigma$. Applying Theorem~\ref{t:stand-transl}
%once again, we conclude that $\mathfrak{M}',w'\models \varphi$
%for all $\varphi\in\Sigma$.
%\endproof
%
%\item[Con] $\ \ \forall x\forall U.(x\in U\to\forall y\forall V.(y\in
%V\to (\exists z.(z\in U\wedge z\in V)\vee \exists z'.(\neg(z'\in
%U)\wedge \neg(z'\in V)))))$

%}}}

%{{{ Van-Benthem characterization

\section{The basic modal language}\label{s:definability}

The expressive power of the basic modal language \emph{on relational
  structures} is relatively well understood. The Van Benthem theorem
characterizes the modally definable \emph{properties of points in
  Kripke models}, in terms of bisimulations, while the
Goldblatt-Thomason theorem characterizes modal definability of
\emph{classes of Kripke frames}, in terms of closure under operations
such as disjoint union.

In this section, we will prove topological analogs of these
results. First, we present a topological version of Van Benthem's
theorem, using the notion of \emph{topo-bisimulations}
\cite{aiello02}.  Next, we identify four operations on topological
spaces that preserve validity of modal formulas. Finally, we
apply these closure conditions in order to determine which
$\mathcal{L}_t$-definable classes are modally definable, and vice
versa.

\subsection{Topological bisimulations}\label{s:topo-bisimulations}

In this section we characterize the modal fragment of $\mathcal L_t$
in terms of topo-bisimulations.

\begin{defn} \label{d:topo-bisim-johan}
  Consider topological models $\mathfrak M=(X,\nu)$ and $\mathfrak
  M'=(X', \nu')$. A non-empty relation $Z\subseteq X\times X'$ is a
  {\it topo-bisimulation between $\mathfrak M$ and $\mathfrak M'$} if
  the following conditions are met for all $x\in X$ and $x'\in X'$:

\begin{description}
\item[\textbf{Zig}] If $x Z x'$ and $x\in O\in\tau$ then there exists $O'\in\tau'$
such that $x'\in O'$ and for all $y'\in O'$ there exists a $y\in O$
such that $yZ y'$.
\item[\textbf{Zag}] If $x Z x'$ and $x'\in O'\in\tau'$ then there exists $O\in\tau$
such that $x\in O$ and for all $y\in O$ there is a $y'\in O'$ such
that $yZy'$.
\item[\textbf{Atom}] If $x Z x'$ then $x \in \nu(p)$ iff $x' \in
\nu'(p)$ for all $p \in \textsc{Prop}$.
\end{description}
Elements $x\in X$ and $x'\in X'$ are said to be \emph{bisimilar},
denoted by $(\mathfrak M,x)\bisim (\mathfrak M',x')$, if there
exists a bisimulation $Z$ between $\mathfrak M$ and $\mathfrak M'$
such that $xZx'$.
\end{defn}

This definition can be formulated more naturally if we use some
standard mathematical notation.
For a binary relation $Z\subseteq X\times X'$ and a set $A\subseteq
X$, let us denote by $Z[A]$ the \emph{image} $\set{x'\in X' \mid
  \exists x\in A.(xZx')}$, and let us define the \emph{preimage}
$Z^{-1}[A']$ of a set $A'\subseteq X'$ analogously.

\begin{prop}\label{p:topo-bis-alt}
  The \textbf{Zig} and \textbf{Zag} conditions in
  Definition~\ref{d:topo-bisim-johan} are equivalent to the following:
\begin{description}
\item[\textbf{Zig$'$}] For all $O\in\tau$, $Z[O]\in\tau'$.
\item[\textbf{Zag$'$}] For all $O'\in\tau$, $Z^{-1}[O']\in\tau$.
\end{description}
\end{prop}
\proof We will only show the equivalence for \textbf{Zig$'$}, the
proof for \textbf{Zag$'$} is analogous. In one direction, suppose
that $Z$ satisfies \textbf{Zig}, and take an open $O \in \tau$.  The
\textbf{Zig} condition ensures that, for each $x'\in Z[O]$, we can
find an open neighborhood $O'\in\tau'$ with $x'\in O'$, such that
$O'\subseteq Z[O]$.  It follows that $Z[O]$, being the union of
these neighborhoods, is open in $\tau'$. For the other direction,
suppose $Z[O]\in\tau'$ holds for all $O\in\tau$.  Consider an
arbitrary $x \in O \in \tau$ and $x' \in X'$ such that $x Z x'$.
Then $Z[O]$ qualifies for an open neighborhood $O'$ of $x'$
satisfying the condition \textbf{Zig} since $x'\in Z[O]\in\tau'$.
\endproof

In what follows, we will freely use this equivalent formulation
whenever it is convenient. Topo-bisimulations are closely linked with the notion of modal
equivalence.

\begin{defn}
  We say that two pointed topological models $(\mathfrak M, x)$ and
  $(\mathfrak M',x')$ are {\it modally equivalent} and write
  $({\mathfrak M}, x)\leftrightsquigarrow ({\mathfrak M'}, x')$ if for
  all formulas $\phi\in\mathcal{ML}$, $({\mathfrak M}, x)\models\phi$
  iff $({\mathfrak M'}, x')\models\phi$.
\end{defn}

\begin{thm}[\cite{aiello02}]\label{t:bisim}
For arbitrary topological pointed models $({\mathfrak M}, x)$ and
$({\mathfrak M'}, x')$, if
$({\mathfrak M}, x)\bisim ({\mathfrak M'}, x')$ then
$({\mathfrak M}, x)\leftrightsquigarrow ({\mathfrak M'}, x')$.
\end{thm}
\proof The proof proceeds via straightforward induction on the
complexity of modal formulas. We only treat the case
$\phi=\Box\psi$.

Suppose $({\mathfrak M}, x)\models\Box\psi$. Then there exists an
open neighborhood $O$ of $x$ such that $O\models\psi$. By
$\textbf{Zig}$ we obtain that $Z[O]$ is an open neighborhood of $x'$
and, by induction hypothesis, $Z[O]\models\psi$. Therefore
$({\mathfrak M'}, x')\models\Box\psi$. The other direction is proved
similarly.
\endproof

The converse does not hold in general, but it holds on a restricted
class of $\mathcal{L}_t$-saturated topological models.

\begin{thm}\label{t:henn-milner}
  Let $\mathfrak M$ and $\mathfrak M'$ be
  $\mathcal{L}_t$-saturated topological models, and suppose
  that $(\mathfrak{M},x)\leftrightsquigarrow (\mathfrak{M}',x')$.
  Then $(\mathfrak{M},x)\bisim (\mathfrak{M}',x')$.
\end{thm}

\proof Let $\mathfrak{M}=(X,\tau,\nu)$ and $\mathfrak{M}' =
(X',\tau',\nu')$, and let $Z\subseteq X\times X'$ be the modal
indistinguishability relation (i.e., $xZx'$ iff $(\mathfrak M,
x)\leftrightsquigarrow (\mathfrak M', x')$).  We will show that $Z$
is a topo-bisimulation, and hence, $(\mathfrak{M},x)\bisim
(\mathfrak{M}',x')$. That the \textbf{Atom} condition holds follows
immediately from the construction of $Z$. In the remainder of this
proof, we will show that $\textbf{Zag}$ holds. The case for
\textbf{Zig} is analogous.

Consider any $a,a'$ such that $aZa'$, and let $O'\in\tau'$ be an
open neighborhood of $a'$. Since $\mathfrak{M}'$ is
$\mathcal{L}_t$-saturated, we may assume that $O'$ is
point-saturated (if not, just take a point-saturated subneighborhood
of $a'$).  We need to find an open neighborhood $O$ of $a$ such that
for each $b\in O$ there exists a $b'\in O'$ with $bZb'$.

By $\mathcal{L}_t$-saturatedness of $\mathfrak{M}$, we know that $a$
has an open neighborhood $O_a$ such that, for every modal formula
$\phi$, if $a\models\Box\phi$ then $\phi$ holds throughout $O_a$.
Dually, this means that

\begin{center}
  (*) \qquad For any $b\in O_a$ and modal formula $\phi$, if
  $b\models\phi$ then $a\models\Diamond\phi$.
\end{center}

%\cproof \ Suppose $b\models\phi$ for some $b\in O$, but
%$a\not\models\Diamond\phi$. Then $a\models\Box\neg\phi$ and by
%$b\in O\subseteq O_{\neg\phi}\models\neg\phi$ we obtain that
%$b\models\neg\phi$---a contradiction.\endproof

To show that $O_a$ meets the requirements of the $\textbf{Zag}$
condition, consider any $b\in O_a$. We will find a $b'\in O'$ such
that $bZb'$. Let $\Sigma_b$ be the set of modal formulas true at
$b$. Every finite subset of $\Sigma_b$ is satisfied somewhere in
$O'$. For, consider any finite $\Sigma'\subseteq\Sigma_b$.  Then by
(*), $\mathfrak{M},a\models\Diamond\bigwedge\Sigma'$, and hence
$\mathfrak{M}',a'\models\Diamond\bigwedge\Sigma'$. Therefore
$\bigwedge\Sigma'$ must be satisfied somewhere in $O'$. Recall that
$O'$ is point-saturated. We conclude that there is a point $b'\in
O'$ satisfying $\Sigma_b$. It follows that $(\mathfrak
M,b)\leftrightsquigarrow (\mathfrak M', b')$, and hence $bZb'$.
\endproof

Combining this with Theorem~\ref{thm:Lt-saturation}, we obtain

\begin{thm}\label{t:vB charcterization}
An $\mathcal L_t$-formula $\alpha(x)$ is invariant
under topo-bisimulations iff it is equivalent to the standard
translation of a modal formula.
\end{thm}
\proof Easily adapted from the proof of the van Benthem
Characterization Theorem for relational semantics (see e.g.
\cite[Theorem 2.68]{blackburn01} for details).
\endproof

%}}}

%{{{ Validity-preserving operatins

\subsection{Validity preserving operations}

In this section, we use topo-bisimulations for showing that three
natural operations on topological spaces (topological sums, open
subspaces and interior maps) preserve  validity of modal
formulas.

\subsubsection{Topological sums}

The \emph{topological sum} (also called \emph{disjoint union},
\emph{direct sum}, or \emph{coproduct}) of a family of disjoint
topological spaces $(X_i,\tau_i)_{i\in I}$, denoted by
$\biguplus_{i\in I}(X_i,\tau_i)$, is the space $(X,\tau)$ with
$X=\bigcup_{i\in I} X_i$ and $\tau=\set{O\subseteq X \mid \forall i\in
  I.(O\cap X_i\in\tau_i)}$. For non-disjoint spaces, the topological
sum is obtained by taking appropriate isomorphic copies.  In the
sequel, when working with topological sums, we will tacitly assume
that the spaces involved are disjoint (cf. \cite[pp.
123-126]{engelking}).

%\begin{defn}[Topological sum]
%  The {\em topological sum} $\biguplus_{i\in I}(X_i,\tau_i)$ of a
%  family of disjoint topological spaces is the space
%  $(X,\tau)$ with $X=\bigcup_{i\in I}X_i$ and
%  topology $\tau=\set{O\subseteq X \mid \forall i\in I.(O\cap X_i\in\tau_i)}$.
%\end{defn}
%

\begin{thm}\label{t:top sum}
  Let $({X}_i)_{i\in I}$ be a family of topological spaces and let
  $\phi$ be a modal formula. Then $\biguplus_{i\in I}X_i\models\phi$
  \textrm{iff} $\forall i\in I.({X}_i\models\phi)$
\end{thm}

\proof There is a natural topo-bisimulation $Z_i$ between  $X_i$
and $X=\biguplus\limits_{i\in I}X_i$:
$$Z_i=\set{(x,x)\mid x\in X_i}$$
Suppose that $X_i\models\phi$ for each $i\in I$. In order to show
that $X\models\phi$, consider any valuation $\nu$ and point $x$.
Clearly $x$ must belong to $X_j$ for some $j\in I$.  Let $\nu_j$ be
the restriction of $\nu$ to $X_j$, i.e., $\nu_j(p) = \nu(p)\cap X_j$
for each $p\in\textsc{Prop}$. It is easily seen that $Z_j$ is a
topo-bisimulation between $((X,\nu),x)$ and $((X_j,\nu_j),x)$. Since
$X_j\models\phi$, we obtain by Theorem~\ref{t:bisim} that
$(X,\nu),x\models\phi$. This argument was independent of $\nu$ and
$x$, and therefore we may conclude that $X\models\phi$.

The other direction is established similarly, and follows also
from Theorem~\ref{t:open subspaces} below.
\endproof

The above lemma can immediately be put to use to show that
compactness and connectedness are not modally definable. Recall that
a space is said to be \emph{compact} if any open cover of the space
contains a finite subcover, and a space is said to be \emph{connected}
if it does not contain a proper non-empty subset that is both closed
and open.

\begin{cor}\label{c:conn not def}
The class of connected spaces and the class of compact
spaces are not modally definable.
\end{cor}
\proof In view of Theorem~\ref{t:top sum} it suffices to note that
while each space $X_i=(\set{i},\set{X_i,\emptyset})$ (a singleton
set equipped with the only possible topology) is both connected and
compact, the topological sum ${X}=\biguplus_{i\in \omega}{X}_i$
is neither connected nor compact.
\endproof

Incidentally, the class of connected spaces is definable in a modal
language with the \emph{global modality} \cite{shehtman}. We discuss
the global modality and the connectedness axiom in
Section~\ref{s:global} below.

Typical examples of properties that \emph{are} preserved under taking
disjoint union are \emph{disconnectedness}, as well as $T_0$, $T_1$,
$T_2$, and discreteness.

\subsubsection{Open subspaces}\label{sec:open subspaces}

Given a topological space $(X,\tau)$ and an open subset $O\in\tau$,
there is a natural topology on $O$ induced by $\tau$, namely
$\tau_O=\set{A\subseteq O \mid A\in\tau}$, or, equivalently,
$\tau_O=\set{A\cap O \mid A\in\tau}$ (cf. \cite[pp.
111-112]{engelking}).  An \emph{open subspace} of $X$ is any space
$(O,\tau_O)$ for $O\in\tau,\ O\neq\emptyset$.

\begin{thm}\label{t:open subspaces}
  Let $(X,\tau)$ be a space and $(O,\tau_O)$ an open subspace, and let
  $\phi$ a modal formula. If $X\models\phi$ then $O\models\phi$.
\end{thm}

\proof Suppose $(O,\nu),x\not\models\phi$, for some valuation $\nu$
and point $x\in O$. We can view $\nu$ also as a valuation for $X$.
The inclusion map is then a topo-bisimulation between $((O,\nu),x)$
and $((X,\nu),x)$.  It follows by Theorem~\ref{t:bisim} that
$(X,\nu),x\not\models\phi$.
\endproof

Theorem~\ref{t:open subspaces} provides us with another way to prove
that connectedness is not modally definable: the real line
$\mathbb{R}$ with the usual topology is connected, but its open
subspace $\mathbb{R}\setminus\{0\}$ is not. Using
Theorem~\ref{t:open
  subspaces} we can also show that \emph{disconnectedness} is not
modally definable. We call a space disconnected if it is not
connected. Since the two-point discrete space is disconnected, while
its one-point open subspaces are connected, we obtain the following
corollary:

\begin{cor}
  The class of disconnected spaces is not modally definable.
\end{cor}

Typical examples of properties that \emph{are} preserved under taking
open subspaces are $T_0$, $T_1$, $T_2$, \emph{density-in-itself} and
\emph{being a Baire space}.

\subsubsection{Images of interior maps}

The third operation that we will consider is taking images of
\emph{interior maps} (also known as \emph{continuous open maps}).  A
map $f:{X}_1\to{X}_2$ between topological spaces $(X_1,\tau_1)$ and
$(X_2,\tau_2)$ is said to be {\em open} if $f(O)\in \tau_2$ for each
$O\in\tau_1$ (i.e. images of opens are open), and \emph{continuous}
if $f^{-1}(O)\in \tau_1$ for each $O\in\tau_2$ (i.e. preimages of
opens are open). If $f$ is both open and continuous, it is called an
\emph{interior map}. Note that homeomorphisms are simply bijective
interior maps (cf. \cite[pp. 57-67]{engelking}).

\begin{thm}\label{t:interior images}
  Let $X_1$ and $X_2$ be topological spaces and $f:X_1\to X_2$ a
  surjective interior map. For all modal formulas $\phi$,
  if $X_1\models\phi$ then $X_2\models\phi$.
\end{thm}
\proof By contraposition: suppose $(X_2,\nu_2),x_2\not\models\phi$ for some $\nu_2, x_2$.
 Let $x_1$ be any element of $X_1$ such that $f(x_1)=x_2$ (recall that $f$
 is surjective), and let $\nu_1$ be the valuation on $X_1$ defined by
  $\nu_1(p)=f^{-1}[\nu_1(p)]$. By construction,
  the graph of $f$ is a topo-bisimulation between $((X_1,\nu_1),x_1)$ and
  $((X_2,\nu_2),x_2)$ (cf. Proposition~\ref{p:topo-bis-alt}).
  It follows by Theorem~\ref{t:bisim} that $(X_1,\nu_1),x_1\not\models\phi$.
\endproof

Not many properties of spaces are preserved under taking images of
interior maps.

It is known that the real line $\mathbb R$ with its usual topology
obeys all separation axioms $T_i$ for
$i\in\set{0,D,1,2,3,3{\frac{1}{2}}, 4,5}$. As a corollary of
Theorem~\ref{t:interior images}, we obtain that none of these are
definable in the basic modal language.

\begin{cor}\label{c:t1-t6}
  The separation axioms $T_i$ with $i\in\set{0,D,1,2,3,3{\frac{1}{2}},
    4,5}$ are not definable in the basic modal language.
\end{cor}
\proof Consider the interior map from the real line $\mathbb R$ with
the standard topology onto $X=\set{1,2}$ equipped with the trivial
topology $\tau=\set{\emptyset, X}$, sending the rationals to $1$ and
the irrationals to $2$. It is easy to verify that the reals obey all
separation axioms, while $X$ obeys none. As surjective interior maps
preserve modal validity, none of the separation axioms can be
defined by a formula in the basic modal language.  \endproof

We will show in Section~\ref{s:extensions} that extending the modal
language can help us in defining some of the lower separation
axioms.

Examples of properties that \emph{are} preserved under interior maps
are being HI, extremally disconnected, compact, connected or
separable (in fact, the latter three are even preserved by
continuous maps).

For further application of the preservation results presented in
this section, as well as related techniques for establishing
(un)definability of topological properties such as submaximality,
being nodec, door, maximal, perfectly disconnected, etc., see the
recent paper \cite{submax}.

%}}}

%{{{ "Johan's lemma"

\subsection{Alexandroff extensions}\label{s:alexandroff ext}

In this section, we introduce a fourth operation on topological
spaces---formation of \emph{Alexandroff extensions}. It allows one to
turn arbitrary spaces into Alexandroff spaces. We will show that this
construction reflects the validity of modal formulas, and we will
identify a connection between Alexandroff extensions and topological
ultraproducts.

\begin{defn}[Alexandroff extensions]
  A filter $\mathcal F \subseteq \wp(X)$ over a topological space $(X,\tau)$ is
  called \emph{open} if for all $A \in \mathcal F$, also $\I A \in \mathcal F$.
  The \emph{Alexandroff extension} of a space $(X,\tau)$ is the space
  $X^*=(\uf X,\tau^*)$, where $\uf X$ is the set of ultrafilters
  over $X$, and $\tau^*$ is the topology over $\uf X$ generated by
  the sets of the form $\{\mathfrak u\in \uf X \mid \mathcal F\subseteq \mathfrak u\}$ for $\mathcal F$ an open
  filter over $X$.
\end{defn}

\begin{thm}
  For any space $X$, $X^*$ is Alexandroff.
\end{thm}

\begin{proof} For any point $\mathfrak u\in \Ae{X}$ consider a filter
  $\mathcal F$ generated by all open sets that belong to $\mathfrak
  u$. Then the set $\{\mathfrak{v} \in \Ae{X} \mid \mathcal F
  \subseteq \mathfrak{v}\}$ is a least open neighborhood of
  $\mathfrak u$. It follows that $\mathfrak v$ is in the least open
  neighborhood of $\mathfrak u$ iff for each $\I A\in\mathfrak u$ we
  have $\I A\in\mathfrak v$ iff $\C A\in\mathfrak u$ for each $A\in\mathfrak v$.
\end{proof}

Note that the map $\pi:X\to\Ae{X}$ that sends $a\in X$ to the
corresponding principal ultrafilter $\pi_a$ need not be open, or
even continuous \cite[Example 5.13]{BBM06}. Indeed the image
$\pi(X)$, as a subspace of $\Ae{X}$, might not be homeomorphic to
$X$---as soon as $X$ is $T_1$ the subspace $\pi(X)$ is discrete.
Nevertheless, it is worth mentioning that the topology $\tau^*$
preserves the information about the original topology $\tau$ in a
curious way. It is an easy exercise for the reader familiar with
ultrafilter convergence (see, e.g., \cite[pp. 91-93]{engelking})
that $\mathfrak u\in\Ae{X}$ belongs to the least open neighborhood
of the principal ultrafilter $\pi_a$ in $\Ae{X}$ iff $\mathfrak u\to
a$ (i.e. $\mathfrak u$ converges to $a\in X$ according to $\tau$).

Basic open sets of the Alexandroff extension $\Ae{X}$ have a nice
characterisation that follows immediately from their definition. For
any topological space $X$ and subset $A\subseteq X$, let $\s{A} =
\{\mathfrak u \in X^* \mid A \in \mathfrak u \}$. It easily seen
that:
\begin{itemize}
%% \item $\uf \PP{X}$ equipped with the topology generated by $*$-images
%%   of all open sets of $T$ is compact extension of T (see chapter 3.2
%%   in \cite{gabelaia01})
\item[-] $\s{\{a\}}=\set{\pi_a}$;
\item[-] $\s{(A \cap B)} = \s{A} \cap \s{B}$, $\s{(A \cup B)} = \s{A}
\cup \s{B}$;
\item[-] $A^*$ is open iff $A$ is open.
\end{itemize}
Now, the basic open sets of $X^*$ are precisely the sets of the form
$\bigcap_{A\in\mathcal F} \s{A}$ for $\mathcal{F}$ an open filter on
$X$.

We saw in the earlier sections that some topological constructions
preserve modal validity. Now we show that formation of the
Alexandroff extension \emph{anti-preserves} modal validity.

\begin{thm}\label{t:ae reflects validity}
  Let $X$ be a topological space and $X^*$ its Alexandroff extension.
  For all modal formulas $\phi$, if $X^*\models\phi$ then
  $X\models\phi$.
\end{thm}

\proof By contraposition: suppose $X\not\models\phi$. Then there
exists a valuation $\nu$ such that $\nu(\neg\phi)\neq\emptyset$. Let
$\nu^*$ be the valuation on $X^*$ defined by
$\nu^*(p)=\set{\mathfrak u\in X^*\mid \nu(p)\in\mathfrak u}$.  We
will show that
\begin{center}
(*) \qquad For any $\psi\in\mathcal{ML}$ and $\mathfrak u\in X^*$,
$\mathfrak u\models\psi$ iff $\nu(\psi)\in \mathfrak u$.
\end{center}
This gives us the intended result: since $\nu(\neg\phi)\neq\emptyset$,
we can extend $\nu(\neg\phi)$ to an ultrafilter. It follows that
$\nu^*(\neg\phi)\neq\emptyset$, so $X^*\not\models\phi$, as required.

We will prove (*) by induction on the complexity of the formula
$\psi$.  The propositional case is taken care of by the definition of
$\nu^*$, the cases for the boolean connectives are rather obvious, so
we only address the modality case. Let $\psi$ therefore be of the form
$\Box\xi$.

[$\Rightarrow$] Suppose $\mathfrak u\models\Box\xi$. Then $\mathfrak
u$ has an open neighborhood (restrict to the element of the base
without loss of generality) $O=\{\mathfrak v\in \Ae{X}\mid \mathcal
F\subseteq \mathfrak v\}$ such that $\mathcal F$ is an open filter
over $X$ and $\mathfrak v\models\xi$ holds for all $\mathfrak v\in
O$. In other words,
$$\mathcal F\subseteq \mathfrak v\ \ \ \Rightarrow\ \ \ \ \mathfrak
v\models\xi$$
By the induction hypothesis this can be rephrased as
$$\mathcal F\subseteq \mathfrak v\ \ \ \Rightarrow\ \ \ \ \nu(\xi)\in
\mathfrak v$$
for all $\mathfrak v\in X^*$. This indicates that $\nu(\xi)\in
\mathcal F$.  As $\mathcal F$ is an open filter, we obtain
$\I\nu(\xi)\in \mathcal F$. Since $\mathfrak u$ extends $\mathcal F$, it
follows, that $\I\nu(\xi)=\nu(\Box\xi)\in \mathfrak u$.

[$\Leftarrow$] Suppose $\nu(\Box\xi)\in \mathfrak u$. Then
$\I\nu(\xi)\in \mathfrak u$. Consider any ultrafilter $\mathfrak v$
from the least open neighborhood of $\mathfrak u$. Clearly
$\I\nu(\xi)\in \mathfrak v$. By $\I\nu(\xi)\subseteq \nu(\xi)$ we
get $\nu(\xi)\in \mathfrak v$. By the induction hypothesis
$\mathfrak v\in\nu^*(\xi)$. As $\mathfrak v$ was arbitrarily chosen
from the least open neighborhood of $\mathfrak u$, we arrive at
$\mathfrak u\models\Box\xi$.
\endproof

We can immediately conclude that

\begin{cor}
The class of Alexandroff spaces is not modally definable.
\end{cor}
\proof Indeed, suppose a formula $\alpha$ defines the class of
Alexandroff spaces. Take an arbitrary non-Alexandroff space $X$.
Then $X^*$ is Alexandroff, so $X^*\models\alpha$ and by the above
theorem $X\models\alpha$. It follows that $X$ is Alexandroff,
contrary to our assumption.\endproof

The following key theorem (which can be seen as a topological
analogue of \cite[Theorem 3.17]{blackburn01}) connects Alexandroff
extensions to topological ultrapowers.

\begin{thm} \label{t:up_ae}
  For every topological space $X = (X,
  \tau)$ there exists a topological ultrapower $\widehat{\prod_{\mathfrak D} X}$ and a surjective interior map
  $f: \widehat{\prod_{\mathfrak D} X} \to \Ae{X}$. In a picture:
$$
\xymatrix {\widehat{\prod_{\mathfrak D} X} \ar@{-}[d] \ar@{->}[dr]^f &  \\
X & \Ae{X} \\}
$$
\end{thm}

\proof Let us consider an $\mathcal L^2$-based language containing a
unary predicate $P_A$ for each $A \subseteq X$, interpreted naturally
on $X$, i.e., $(P_A)^X=A$. In what follows we will treat $X$ as a
topological model for this (possibly uncountable) language.  By
Theorem~\ref{thm:Lt-saturation}, $X$ has an $\mathcal{L}_t$-saturated
topological ultrapower $\widehat{\prod_{\mathfrak D} X}$. Denote
$Y\equiv \widehat{\prod_{\mathfrak D} X}$. The following
$\mathcal{L}_t$-sentences are clearly true in $X$, and hence, by
Theorem~\ref{t:Los-Lt}, also in $Y$:

\begin{itemize}
\item[(1)] $Y\models\exists x.P_A(x)$ for each non-empty $A\subseteq X$,
\item[(2)] $Y\models\forall x .(P_A(x) \wedge P_B(x) \leftrightarrow
P_{A \cap B}(x))$ for each $A,B\subseteq X$,
\item[(3)] $Y\models\forall x .(\neg P_A(x) \leftrightarrow
P_{-A}(x))$ for each $A\subseteq X$,
\item[(4)] $Y\models\forall x.(P_{\I A}(x)\leftrightarrow\exists U.(x\varepsilon U\wedge\forall y.(y\varepsilon U\to
P_A(y))))$ for each $A\subseteq X$,
\item[(5)] $Y\models\forall x.(P_{\C A}(x)\leftrightarrow\forall U.(x\varepsilon U\to\exists y.(y\varepsilon U\wedge
P_A(y))))$ for each $A\subseteq X$.
\end{itemize}

We define the desired interior map $f: Y \to \Ae{X}$ in the
following way:
$$
f(a)=\{A \subseteq X \mid a \in (P_A)^Y \}
$$

In the remainder of this proof, we will demonstrate that $f$ is
indeed a surjective interior map from $Y$ to $X^*$. First we show
that $f$ is a well-defined onto map.

\begin{itemize}
\item For any $a \in Y$, $f(a)$ is an \emph{ultrafilter} over
$X$.

Recall that an ultrafilter over $X$ is any set $\mathfrak{u}$ of
subsets of $X$ satisfying (\emph{i})
$A\cap B\in\mathfrak{u}$ iff both $A\in\mathfrak{u}$ and $B\in\mathfrak{u}$,
and (\emph{ii}) $A\in\mathfrak{u}$ iff $(X\setminus A)\not\in\mathfrak{u}$.
By (2) and (3) above, $f(a)$ indeed satisfies these properties.

\item   $f$ is \emph{surjective} (i.e., every ultrafilter over $X$ is $f(a)$ for
  some $a\in Y$).

  Take $\mathfrak u \in \Ae{X}$, and let $\Gamma_{\mathfrak{u}}(x)=\{
  P_A(x) \mid A \in \mathfrak u \}$. It follows from (1) and (2) that
  every finite subset of $\Gamma_{\mathfrak{u}}(x)$ is satisfied by
  some point in $Y$.  Since $Y$ is point-saturated, there
  exists $a \in Y$ satisfying $\Gamma_{\mathfrak{u}}(x)$, hence $f(a)
  = \mathfrak u$.
\end{itemize}

Next we show that $f$ is open and continuous. Note that by
Proposition~\ref{p:topo-bis-alt} it suffices to prove that the graph
of $f$ is a topo-bisimulation.

Take arbitrary $a\in Y$ and let $O_a$ be as described in
Definition~\ref{def:Lt-saturation}. Let $O_{\mathfrak u}$ be a least
open neighborhood of $\mathfrak u=f(a)$. We proceed by verifying the
conditions {\bf Zig} and {\bf Zag} for the pair $(a,\mathfrak u)$.

\begin{itemize}

\item {\bf Zig}.  Take arbitrary $O'$ such that $a\in O'$. By $\mathcal{L}_t$-saturatedness of $Y$,
there exists a \emph{point-saturated} $O\subseteq O'$ such that $a\in O$.

  Take arbitrary $\mathfrak v\in O_{\mathfrak u}$. We will find a $b\in O$ such that $\mathfrak
  v=f(b)$. Let
  \[ \Gamma_{\mathfrak{v}}(x) = \{P_A(x) \mid A\in\mathfrak{v}\} \]
  Every finite subset of $\Gamma_{\mathfrak{v}}(x)$ is satisfied
  somewhere in $O$. Indeed, if $P_{A_1}, \ldots,
  P_{A_n}\in\Gamma_{\mathfrak{v}}$, denote $B\equiv \bigcap_i A_i$.
  Then $B\in\mathfrak v$ and hence $\C B\in\mathfrak u$. Therefore
  $Y\models P_{\C B}(a)$. It follows by (5) that $P_B$ holds somewhere in $O$.
  By the point-saturatedness of $O$ we may conclude that
  some $b\in O$ satisfies all of $\Gamma_{\mathfrak{v}}(x)$,
  and hence $f(b)=\mathfrak{v}$.

\item {\bf Zag}. It suffices to show that for any $b\in O_a$ we have
$f(b)\in O_{\mathfrak u}$. Suppose the contrary. Then we have $b\in
O_a$ and $f(b)\not\in O_{\mathfrak u}$. The latter means that there
exists a set $A\subseteq X$ such that $A\in f(b)$ but $\C A\not\in
\mathfrak u$. From $A\in f(b)$ we obtain $Y\models P_A(b)$. While
$\C A\not\in \mathfrak u$ ~iff~ $-\C A\in \mathfrak u$ ~iff~ $\I
-A\in\mathfrak u$ ~iff~ $Y\models P_{\I -A}(a)$ ~iff~
$Y\models\exists U.[a\varepsilon U\wedge\forall y.(y\varepsilon U\to
P_{-A}(y))]$ ~iff~ $P_{-A}(x)$ is true throughout some open
neighborhood of $a$ ~iff~ $P_{-A}(x)$ is true throughout $O_a$,
which contradicts $Y\models P_A(b)$ since $b\in O_a$.
\end{itemize}

\vspace{-10mm}
\endproof

%}}}

%{{{ Gb-Th model-theretic proof for BML

\subsection{Modal definability vs $\mathcal L_t$-definability}\label{s: GbTh main}

In this section we are seeking \emph{necessary and sufficient}
conditions, in the spirit of the Goldblatt-Thomason theorem, for a
class of topological spaces to be modally definable. We have already
found some necessary conditions: we have seen that every modally
definable class of topological spaces is closed under the formation of
\emph{topological sums}, \emph{open subspaces} and \emph{interior
  images} and reflects \emph{Alexandroff extensions}.  Our aim is to
prove a converse, in other words, to characterize modal definability
in terms of these closure properties.

\begin{thm}\label{t:GbTh basic Lt}
  Let $\mathsf{K}$ be any $\mathcal{L}_t$-definable class of
  topological spaces. Then $\mathsf K$ is modally definable iff it
  is closed under taking open subspaces, interior images, topological
  sums and it reflects Alexandroff extensions.
\end{thm}

\proof We will only prove the difficult right-to-left direction. The
left-to-right direction already follows from theorems \ref{t:top sum},
\ref{t:open subspaces}, \ref{t:interior images} and \ref{t:ae reflects validity}.

Let $\mathsf{K}$ be any class satisfying the given closure conditions.
Take the set $Log(\mathsf{K})$ of modal formulas valid on
$\mathsf{K}$. We will show that, whenever $X\models Log(\mathsf{K})$,
then $X\in\mathsf{K}$. In other words, $Log(\mathsf{K})$ \emph{defines}
$\mathsf{K}$.

Suppose $X\models Log(\mathsf{K})$ for some topological space $X$.
Introduce a propositional letter $p_A$ for each subset $A\subseteq
X$, and let $\nu$ be the natural valuation on $X$ for this (possibly
uncountable) language, i.e. $\nu(p_A)=A$ for all $A\subseteq X$. Let
$\Delta$ be the set of all modal formulas of the following forms
(where $A,B$ range over subsets of $X$):
\[\begin{array}{lll}
  p_{A\cap B}  &\leftrightarrow& p_A\land p_B \\
  p_{-A}      &\leftrightarrow& \neg p_A \\
  p_{\I A}    &\leftrightarrow& \Box p_A \\
  p_{\C A}    &\leftrightarrow& \Diamond p_A
\end{array}\]
By definition, $\Delta$ is valid on $\mathfrak M=(X,\nu)$. Note that
the standard translations of the formulas in $\Delta$ correspond
exactly to the formulas listed in conditions (2)--(5) of
Theorem~\ref{t:up_ae} (in the corresponding $\mathcal{L}_t$-language,
which has a one-place predicate $P_A(x)$ for each $A\subseteq X$).
What is missing is the condition (1). The following claim addresses
this.

\begin{quote}
\emph{Claim:} For each $a\in X$ there is a model
  $\mathfrak{N}_a=(Y_a,\mu_a)$ with $Y_a\in\mathsf{K}$, such that
  $\mathfrak{N}_a\models\Delta$ and some point in $\mathfrak{N}_a$
  satisfies $p_a$.

  \medskip \cproof Take any $a\in X$, and let
  $\Delta_a=\set{\Box\varphi\mid\varphi\in\Delta} \cup
  \set{p_{\set{a}}}$. As a first step, we will show that there is a
  topological model $\mathfrak{K}$ based on a space in $\mathsf{K}$,
  such that some point $a'$ of $\mathfrak{K}$ satisfies $\Delta_a$. By
  the compactness of $\mathcal{L}_t$ (Theorem~\ref{thm:compactness}),
  it suffices to show that every finite conjunction $\delta$ of
  formulas in $\Delta_a$ is satisfiable on $\mathsf{K}$. Since
  $\delta$ is satisfied at $a$ in $\mathfrak M$ and $\mathfrak
  M\models Log(\mathsf{K})$, $\neg\delta$ cannot belong to
  $Log(\mathsf{K})$. Hence $\delta$ is satisfiable on $\mathsf{K}$.

  By Theorem~\ref{thm:Lt-saturation} we may assume $\mathfrak{K}$ is
  $\mathcal{L}_t$-saturated. Let $O_{a'}$ be an open neighborhood of
  $a'$ as described in Definition~\ref{def:Lt-saturation}, and let
  $\mathfrak{N}_a$ be the submodel of $\mathfrak K$ based on $O_{a'}$.
  Then $\mathfrak{N}_a$ satisfies all requirements of the claim.
  \endproof
\end{quote}

Note how, in the above argument, we used the fact that $\mathsf{K}$ is
$\mathcal{L}_t$-definable (for the compactness argument, and for the
saturation), and that it is closed under taking open subspaces.  Next,
we will use the fact that $\mathsf{K}$ is closed under taking
\emph{topological sums}.

Let $Y = \biguplus_{a\in X} Y_a$, and let $\mathfrak{N}=(Y,\mu)$,
where $\mu$ is obtained from the $\mu_a$'s in the obvious way. By
closure under taking topological sums, $Y\in\mathsf{K}$. Moreover, by
Theorem~\ref{t:top sum}, $\mathfrak{N}\models\Delta$. Finally, each
$p_A$, for non-empty $A\subseteq X$, holds at some point in
$\mathfrak{N}$ (more precisely, at some point in $\mathfrak{N}_a$ for
any $a\in A$).  It follows (using the standard translation) that the
conditions (1)--(5) from the proof of Theorem~\ref{t:up_ae} hold for
$\mathfrak N$.

We can now proceed as in the proof of Theorem~\ref{t:up_ae}, and
construct an interior map from an ultrapower of $\mathfrak N$ onto the
Alexandroff extension $X^*$ of $X$.  Since $\mathsf{K}$ is closed
under topological ultrapowers (Theorem~\ref{t:Los-Lt}) and images of
interior maps, and reflecs Alexandroff extensions, we conclude that
$X\in \mathsf{K}$.
\endproof

Inspection of the proof shows that Theorem~\ref{t:GbTh basic Lt}
applies not only to $\mathcal L_t$-definable classes but to any class
of spaces closed under ultraproducts.  In fact, by Lemma~\ref{l:u pr
  -> u power top sum} below, closure under ultra\emph{powers}
suffices. Some further improvements are still possible. Most
importantly, using algebraic techniques, we will show in the next
section that \emph{closure under Alexandroff extensions} already
suffices. For the complete picture, see Corollary~\ref{c:GbTh basic full}.
%
% (we could have taken the
% ultraproduct of $Y_a$ instead of their topological sum in the proof).
%
% Why is this relevant????
%
%Since ultraproducts are interior images of box products, it also
%applies to any class closed under box products.

\paragraph{The opposite question.}

Theorem~\ref{t:GbTh basic Lt} characterizes, among all
$\mathcal{L}_t$-definable classes of topological spaces, those that
are modally definable.  It makes sense to ask the opposite question:
\emph{which modally definable classes of spaces are
  $\mathcal{L}_t$-definable?} In classical modal logic the answer was
provided by van Benthem in \cite{johan76} (see also \cite{rob75}).  We
follow the route paved in these papers. First we prove a topological
analogue of an observation due to Goldblatt:

\begin{lem}\label{l:u pr -> u power top sum}
An ultraproduct of topological spaces is homeomorphic to an open
subspace of the ultrapower (over the same ultrafilter) of their
topological sum.
\end{lem}
\proof Suppose $(X_i)_{i\in I}$ is a family of topological spaces
and $\mathfrak D$ is an ultrafilter over $I$. Denote by
$X=\biguplus_{i\in I} X_i$ the topological sum of $X_i$ and by
$Y=\widehat{\prod_{\mathfrak D} X_i}$ their topological
ultraproduct. Take arbitrary $a: I\to \biguplus_{i\in I} X_i$ such
that $a(i)\in X_i$. Then $a$ can be viewed both as an element of
$\prod_{i\in I} X_i$ and as an element of $\prod_{i\in I} X$. This
defines a natural embedding from $Y$ into $\widehat{\prod_{\mathfrak
D} X}$ which is clearly injective. That this embedding is open is
easily seen (recall that this suffices to be checked on the elements
of the base). To show that it is also continuous, suppose
$[a]_\mathfrak D\in\prod_{\mathfrak D} X$ is such that $A=\set{i\mid
a(i)\in X_i}\in\mathfrak D$ (so $[a]_\mathfrak D$ comes from $Y$).
Then any basic ultrabox neighborhood $\prod_{\mathfrak D} O_i$ of
$[a]_\mathfrak D$ is such that $B=\set{i\mid a(i)\in O_i\subseteq
X}\in\mathfrak D$. We clearly have $A\cap B\in \mathfrak D$, so
$\prod_{\mathfrak D} (O_i\cap X_i)$ is another open neighborhood of
$[a]_\mathfrak D$, now also in $Y$. The required continuity follows.
Since we have established that $Y$ can be embedded into
$\widehat{\prod_{\mathfrak D} X}$ by an interior map, it follows
that $Y$ is homeomorphic to an open subspace of
$\widehat{\prod_{\mathfrak D} X}$.
\endproof

We are one step away from finding a nice criterion for a modally
definable class to be $\mathcal{L}_t$-definable. It follows from
Garavaglia's theorem \cite{garavaglia} (the topological analogue of
the Keisler-Shelah Theorem) that a class $\mathsf K$ of spaces is
$\mathcal{L}_t$-definable iff $\mathsf K$ is closed under
isomorphisms and ultraproducts and the complement of $\mathsf K$ is
closed under ultrapowers.

\begin{thm}
A modally definable class $\mathsf K$ of spaces is $\mathcal{L}_t$-definable iff it
is closed under ultrapowers.
\end{thm}
\proof If $\mathsf K$ is $\mathcal{L}_t$-definable, then it is
clearly closed under ultrapowers. For the converse direction take a
modally definable class $\mathsf K$ that is closed under
ultrapowers.  Then $\mathsf K$ is closed under topological sums and
open subspaces. It follows from Lemma~\ref{l:u pr -> u power top
sum} that $\mathsf K$ is closed under ultraproducts. It is easily
seen that any modally definable class is closed under $\mathcal
L_t$-isomorphisms. It follows from
Theorems~\ref{t:Los-Lt}~and~\ref{t:stand-transl} that the complement
of $\mathsf K$ is closed under ultrapowers. Hence $\mathsf K$ is
$\mathcal{L}_t$-definable.\endproof

Since modally definable classes are closed under interior images and
ultrapowers are interior images of box products via the canonical
quotient map, we obtain

\begin{cor}
A modally definable class of spaces that is closed under box powers
is $\mathcal{L}_t$-definable.
\end{cor}

\paragraph{Separating examples.}
To close this section we give examples separating
$\mathcal{L}_t$-definability from modal definability.  We have
exhibited earlier $\mathcal L_t$-sentences defining the separation
axioms $T_0-T_2$. We have also shown in Corollary~\ref{c:t1-t6} that
$T_0-T_2$ are not definable in the basic modal language. Thus we
have examples of $\mathcal{L}_t$-definable classes of spaces that
are not modally definable.

To show that there are modally definable classes of spaces that are
not $\mathcal{L}_t$-definable requires more work. Recall that the
class of Hereditarily Irresolvable (HI) spaces is modally definable
(Theorem~\ref{t:HI definable}). This class is not $\mathcal
L_t$-definable. We will use the following lemma:

\begin{lem} Any class $\mathsf{K}$ of spaces that is both modally
  definable and $\mathcal L_t$-definable is closed under Alexandroff
  extensions.
\end{lem}
\proof Suppose $X\in \mathsf{K}$. By Theorem~\ref{t:up_ae} there
exists a topological ultrapower $Y$ of $X$ and an onto interior map
$f:Y\to \Ae{X}$. Being an $\mathcal L_t$-definable class, $\mathsf{K}$
is closed under topological ultrapowers.  Hence
$Y\in\mathsf{K}$. Being a modally definable class, $\mathsf{K}$ is
closed under interior images. Therefore $\Ae{X}\in\mathsf K$, as
required.\endproof

\begin{thm}The class of HI spaces is not $\mathcal{L}_t$-definable.
\end{thm}

\proof By Theorem~\ref{t:HI definable} and the above lemma, to prove
that the class of HI spaces is not $\mathcal L_t$-definable it
suffices to show that this class is not closed under Alexandroff
extensions. In \cite[Example 5.12]{BBM06} a space $X$ is exhibited
that is HI, but its Alexandroff extension is not HI. We reproduce
this example for reader's convenience.

Let $X=(\N, \tau)$ be a topological space with carrier $\N=\set{1,
2,\dots}$ and topology $\tau=\set{[1,n)\mid n\in\N}\cup\N$. This is
the Alexandroff topology corresponding to the order $\geq$. To show that $X$ is HI,
observe first that for an arbitrary subset $A\subseteq \N$ we have $\C A=[\min
A,\infty)$. Further, if $A,A'\subseteq B$ are such that $A\cap
A'=\emptyset$ it is easily seen that either $\min A> \min B$ or $\min A'> \min
B$. Hence either $B\not\subseteq \C A$ or $B\not\subseteq \C A'$.
This shows that no subset of $X$ can be decomposed into two disjoint dense in
it sets, so $X$ is HI.

Consider the Alexandroff extension $\Ae{X}$. Let $\mathfrak
F\subseteq \Ae{X}$ denote the set of all the free ultrafilters over
$X$. Fix two distinct free ultrafilters $\mathfrak u,\mathfrak
v\in\mathfrak F$. We will show that both $\mathfrak u$ and
$\mathfrak v$ belong to the least open neighborhood of any
$\mathfrak w\in\mathfrak F$. To see this it suffices to check that
for any non-empty $A\subseteq X$ we have $\C A\in\mathfrak w$. But
if $A$ is non-empty, then $\C A=[\min A,\infty)$ is cofinite and
thus belongs to the free ultrafilter $\mathfrak w$. It follows that
$\set{\mathfrak u}$ and $\set{\mathfrak v}$ are two
 disjoint dense in $\mathfrak F$ subsets. Hence $\Ae{X}$ is not HI.
\endproof

%\todo{Is it $\mathcal{L}^2$ definable? Perhaps put as a question.  It
%  is natural to ask whether there are modally definable classes that
%  are not $\mathcal{L}^2$-definable. In $\mathcal{L}^2$ one can only
%  quantify over open sets, while in modal logic, validity is defined
%  via full powerset quantification. Would ``weakly hereditarily
%  irresolvable'' be a weaker property?}

\section{Interlude: an algebraic perspective}\label{s:algebra}

In this paper we have chosen to approach the question of
definability from the \emph{model-theoretic} perspective. While this
approach is rather powerful and fruitful, it is not the only
possible one. In this section we sketch an equally potent approach
via \emph{Universal Algebra}. We will outline how, using algebraic
techniques, one can prove a slightly stronger version of
Theorem~\ref{t:GbTh basic Lt}. It should be noted however, that the
algebraic techniques do not straightforwardly generalize to various
extensions of the basic modal language. The model theoretic approach
provides more flexibility in this respect, as we will see in
Section~\ref{s:extensions}.

Most of the proofs that are missing in this section can be
found in \cite{gabelaia01}.

\subsection{Algebraic semantics for modal logics}

In a certain sense, the algebraic semantics for modal logic is most
adequate, however it is also most abstract.
Here we give a basic intuition of the universal algebraic approach
to modal logic. More details can be found in standard textbooks
\cite{cz, blackburn01}

\begin{defn}\label{d:modal algebra}
A \emph{modal algebra} is a tuple $(B, \Box)$ where $B$ is a Boolean
algebra and $\Box: B\to B$ is an operator such that for all $a,b\in
B$ the following holds:

\vspace{2mm}

\begin{tabular}{ll}
(i1) & $\Box \top=\top$,\\
(i2) & $\Box(a\wedge b)= \Box a\wedge \Box b$.\\
\end{tabular}
\end{defn}

It is easily seen how $\mathcal{ML}$ can be interpreted on a modal
algebra. Propositional letters designate elements of the Boolean
algebra and the operations are interpreted by their algebraic
counterparts. Every modal formula then becomes a polynomial that can
be computed on tuples of elements of the algebra. The formulas that
evaluate to $\top$ regardless of the assignment of the elements to
the propositional letters are said to be \emph{valid} in the
algebra. It can be shown that every class of modal algebras
determines a normal modal logic of the formulas that are valid on
every algebra of the class \cite[Chapter 7]{cz}. Conversely, a
normal modal logic singles out a class (indeed, a \emph{variety}) of
modal algebras that validate all the formulas in the logic. This
correspondence between varieties and logics is 1-1 \cite[Chapter
7]{cz}.

Modal algebras arising from topological spaces are called
interior algebras. We discuss them next.

\subsection{Interior algebras}

Each non-empty set $X$ gives rise to the Boolean algebra $\wp{X}$ of
its subsets. Suppose in addition $X$ is endowed with a topology $\tau$. How is it possible to represent this additional
structure algebraically? One natural possibility is to consider the
Boolean algebra of all subsets with the corresponding interior operator
$(\wp{X},\I)$. It is known that operation of interior satisfies the well-known
Kuratowski axioms \cite{engelking}. In fact, topological spaces can
equivalently be described as sets endowed with operators satisfying the following:\\

\begin{tabular}{ll}
(I1) & $\I X=X$ \\
(I2) & $\I A\subseteq A$ \\
(I3) & $\I (A\cap B)=\I A\cap \I B$ \\
(I4) & $\I \I A=\I A$\\
\end{tabular}\\

Abstracting away from powerset Boolean algebras brings us to

\begin{defn}\label{d:interior algebra}
An \emph{interior algebra} is a modal algebra $(B, \Box)$ such that for all $a\in
B$ the following holds:

\vspace{2mm}

\begin{tabular}{ll}
(i3) & $\Box a\leq a$,\\
(i4) & $\Box\Box a= \Box a$.\\
\end{tabular}
\end{defn}

Interior algebra homomorphisms are Boolean homomorphisms that
commute with $\Box$.

More details on interior algebras are contained in \cite{rasiowa&sikorski, blok76}.

\subsection{Duality between spaces and interior algebras}

It has been indicated above that each topological space $X$
naturally gives rise to an interior algebra
$$X^+=(\wp{X},\I)$$
We call $X^+$ the \emph{complex algebra} of $X$.

In fact the map $(\cdot)^+$ can be extended to do more than just producing an interior algebra from a topological
space. Given an interior map $f:X\to
Y$ between two topological spaces we can naturally manufacture an
interior algebra homomorphism $f^+:Y^+\to X^+$ by putting
$f^+=f^{-1}$. Furthermore, it can easily be checked that the map
$(\cdot)^+$ takes the topological sum of spaces into the algebraic product of
the corresponding complex algebras. Thus $(\cdot)^+$ witnesses half of a duality going from
topological spaces to interior algebras.

Another direction of the duality is provided by the construction of Alexandroff
extensions of interior algebras. This is a straightforward
generalization of the corresponding construction for topological
spaces.

\begin{defn}[Alexandroff extensions of interior algebars]
  The \emph{Alexandroff extension} of an interior algebra $(B,\Box)$ is the space
  $B_+=(\uf B,\tau_+)$, where $\uf B$ is the set of ultrafilters
  of $B$, and $\tau_+$ is the topology over $\uf B$ generated by
  the sets of the form $\mathcal F^*=\{\mathfrak u\in \uf B \mid \mathcal F\subseteq \mathfrak u\}$ for $\mathcal F$ an open
  filter over $X$.
\end{defn}

Here by an open filter we mean a filter $\mathcal F$ of the Boolean
algebra $B$ such that if $a\in\mathcal F$ then $\Box a \in\mathcal
F$.

Again it can easily be demonstrated that whenever $h:B\to C$ is an
injective (surjective) interior algebra homomorphism, then the surjective (injective) interior map
$h_+:C_+\to B_+$ can naturally be defined by putting $h_+(\mathfrak u)=\set{a\in {B}\mid h(a)\in \mathfrak u}$.

The maps $(\cdot)^+$ and $(\cdot)_+$ provide us with a link
(duality) between interior algebras and homomorphisms on the one
hand and topological spaces and interior maps on the other. With the
help of this duality we can transfer the question \emph{`which
classes of topological spaces are modally definable?'} to the domain
of interior algebras, where it obtains the following form:
\emph{`which classes of interior algebras are equationally
definable?'} and is immediately answered by the fundamental theorem
of Birkhoff---\emph{`those and only those that are closed under
products, subalgebras and homomorphic images'}.

We have just outlined the proof of the following.

\begin{thm}\label{t:GbTh basic ae}
The class $\mathsf K$ of topological spaces which is closed under
the formation of Alexandroff extensions is modally definable iff it
is closed under taking open subspaces, interior images, topological
sums and it reflects Alexandroff extensions.
\end{thm}

The proof of this theorem, as well as the details of the duality
sketched above are presented in \cite{gabelaia01}. Another
characterization of the modal definability for topological spaces
that applies to \emph{any} class of spaces is contained in
\cite[Theorem 5.10]{BBM06} and is also based on the duality outlined
in this section. For the co-algebraic perspective on modal
definability that encompasses both the relational and the
topological cases, as well as more general semantical frameworks, we
refer to \cite{kurz}.

Combining Theorem~\ref{t:GbTh basic ae} with Theorem~\ref{t:up_ae},
we obtain our most general version of the definability theorem for
the basic modal language:

\begin{cor}\label{c:GbTh basic full}
Let $\mathsf{K}$ be a class of topological spaces satisfying at
least one of the following conditions:
\begin{itemize}
\item[(i)] $\mathsf K$ is $\mathcal{L}_t$-definable;
\item[(ii)] $\mathsf K$ is closed under box powers;
\item[(iii)] $\mathsf K$ is closed under ultrapowers;
\item[(iv)]  $\mathsf K$ is closed under Alexandroff extensions.
\end{itemize}
then $\mathsf K$ is modally definable iff it is closed under taking
open subspaces, interior images, topological sums and it reflects
Alexandroff extensions.
\end{cor}
\proof The easier `only if' part follows from theorems \ref{t:top sum},
\ref{t:open subspaces}, \ref{t:interior images} and \ref{t:ae reflects validity}.

To prove the `if' part suppose that
$\mathsf K$ is closed under taking open subspaces, interior images,
topological sums and it reflects Alexandroff extensions. Let us prove
that under these conditions, if $\mathsf K$ satisfies any of the
conditions (i)-(iii) above, then it also satisfies the condition
(iv).

First we show that each of (i) and (ii) implies (iii). Indeed, if
$\mathsf K$ is $\mathcal{L}_t$-definable, then it is closed under
ultrapowers; also, if $\mathsf K$ is closed under box powers, since
ultrapowers are interior images of box powers under the canonical
quotient map and $\mathsf K$ is  closed under interior images, we
obtain that $\mathsf K$ is closed under ultrapowers.

Next we show that (iii) implies (iv). Indeed, it follows from Theorem~\ref{t:up_ae} and the
closure under interior images that if $\mathsf K$ is closed under
taking ultrapowers, then $\mathsf K$ is closed under Alexandroff
extensions.

Thus, in any of the cases (i)-(iv), $\mathsf K$ is closed under Alexandroff extensions.
Now apply Theorem~\ref{t:GbTh basic ae}.
\endproof

The analogue of Theorem~\ref{t:up_ae} for relational semantics has a
neat algebraic proof \cite{GHV}. A similar proof for the
topological case is lacking and we leave this as a challenge for
the interested reader.

%}}}}}}}

%{{{ Gb-Th proof for hybrid languages

\section{Extended modal languages}\label{s:extensions}

In order to increase the topological expressive power of the basic
modal language, various extensions have been proposed. For instance,
Shehtman \cite{shehtman} showed that \emph{connectedness} becomes
definable when the basic modal language is enriched with the
\emph{global modality}. Similarly, $T_0$, $T_1$ and
\emph{density-in-itself} become definable when we enrich the basic
modal language either with nominals or with the difference modality.
In this section, we show \emph{exactly how much definable power} we
gain by these additions, by giving analogues of Theorem~\ref{t:GbTh
  basic Lt} for these extended languages.
Our findings are summarized in Table~\ref{tab:overview-GbTh}
and~\ref{tab:examples} on page~\pageref{tab:overview-GbTh}
and~\pageref{tab:examples}.

We believe Theorem~\ref{t:vB charcterization} could also be
generalized to the languages studied in this section, using
appropriate analogues of topo-bisimulations. However, we have
decided not to pursue this here suspecting the lack of many new
insights.

%Our aim in this section is to delineate the boundaries of the
%expressive power of hybrid modal languages over topological spaces.
%Our characterization of definability for hybrid modal languages is a
%topological reconstruction of the analogous results for the
%relational semantics \cite{cate04}.

\subsection{The global modality}\label{s:global}

In the basic modal language with $\Diamond$ and $\Box$, one can only
make statements about points that are arbitrarily close to the
current point of evaluation. It appears impossible to say, for
instance, that \emph{there is a point satisfying $p$} (i.e., to
express non-emptyness of $p$). The \emph{global modality}, denoted
by $E$, gives us the ability to make such global statements. For
example, $Ep$ expresses non-emptyness of the set $p$, and $A(p\to
q)$ expresses that $p$ is contained in $q$.

Formally, $\mathcal{M}(E)$ extends  the basic
modal language with an extra operator $E$ that has the following
semantics:
$$
\begin{array}{rcl}
\mathfrak{M},w \models E\phi & \mathrm{iff} & \exists v \in X.\
(\mathfrak{M},v \models \phi)\\
\end{array}
$$
The dual of $E$ is denoted by $A$, i.e., $A\varphi$ is short for $\neg
E \neg \varphi$.  The standard translation can be extended in a
straightforward way, by letting $ST_x(E \phi) = \exists
x.(ST_x(\phi))$.  In other words, $\mathcal{M}(E)$ is still a fragment
of $\mathcal{L}_t$.

Shehtman \cite{shehtman} showed that connectedness can be defined
using the global modality:

\begin{prop}
  $A(\Box p \lor \Box \neg p) \to Ap \lor A\neg p$ defines
  connectedness.
\end{prop}

As connectedness is not definable in the basic modal language
(Corollary~\ref{c:conn not def}), this shows that $\mathcal{M}(E)$ is
more expressive than the basic modal language.  As a consequence of
this increased expressive power, certain operations on spaces do not
preserve validity anymore.

\begin{prop}
  Taking open subspaces, or taking topological sums, in general
  does not preserve validity of $\mathcal{M}(E)$-formulas.
\end{prop}

\proof It suffices to show that connectedness is not preserved by
these two operations.  The real interval $(0,1)$, with the usual
topology, is connected, but its open subspace $(0,{\frac{1}{2}})\cup
({\frac{1}{2}}, 1)$ is not. Likewise, for any connected space $X$, the
topological sum $X \biguplus X$ is no longer connected.
\endproof

Taking interior images, on the other hand, \emph{does} preserve
validity of $\mathcal{M}(E)$-formulas, and taking Alexandroff
extensions anti-preserves it. In fact, these two operations
characterize definability in $\mathcal{M}(E)$, as the following
analogue of Theorem~\ref{t:GbTh basic Lt} shows.

\begin{thm}\label{t:GbTh global}
  Let $\mathsf{K}$ be any $\mathcal{L}_t$-definable class of
  topological spaces. Then $\mathsf K$ is definable in the basic modal
  language with global modality iff it is closed under interior images
  and it reflects Alexandroff extensions.
\end{thm}

\begin{proof}
  The `only if' direction is a straightforward adaptation of
  Theorems~\ref{t:interior images} and \ref{t:ae reflects validity}.
  The proof of the `if' direction is essentially a simplification of
  the proof of Theorem~\ref{t:GbTh basic Lt}:
  suppose $X$ is a topological space validating the
  $\mathcal{M}(E)$-theory of $\mathsf{K}$, and let $\Delta$ be the
  following set of formulas, for all $B,C\subseteq X$:

\begin{center}\begin{tabular}{l}
$E p_B$ for non-empty $B$\\
$A(p_{B\cap C} \leftrightarrow p_B\land p_C$)\\
$A(p_{-B} \leftrightarrow \neg p_B)$ \\
$A(p_{\I B} \leftrightarrow \Box p_B)$\\
$A(p_{\C B} \leftrightarrow \Diamond p_B)$
\end{tabular}\end{center}
Note that these formulas exactly correspond to conditions (1)--(5)
from the proof of Theorem~\ref{t:up_ae}.  As in the proof of
Theorem~\ref{t:GbTh basic Lt}, we can find a topological model
$\mathfrak{N}=(Y,\mu)$ with $Y\in\mathsf{K}$, such that
$\mathfrak{N}\models\Delta$. Finally, we proceed as in the proof of
Theorem~\ref{t:up_ae}, and construct an interior map from an
ultrapower of $\mathfrak{N}$ onto the Alexandroff extension $X^*$ of
$X$. Since $\mathsf{K}$ is closed under topological ultrapowers
(Theorem~\ref{t:Los-Lt}) and images of interior maps, and reflects Alexandroff
extensions, we conclude that $X\in\mathsf{K}$.
\end{proof}

\subsection{Nominals}\label{s:nominals}

Another natural extension of the basic modal language is with
\emph{nominals}. Nominals are propositional variables that denote
singleton sets, i.e., they name points.  In point-set topology one
often finds definitions that involve both open sets and individual
points. In the language $\mathcal L_t$, one can refer to the points
in the space by means of point variables. The basic modal language
lacks such means, and nominals can be seen as a way to solve this
problem. Here are some examples of properties that can be defined
using nominals:
$$
\begin{array}{ll}
T_0 & @_i \Diamond j \land @_j\Diamond i \to @_ij \\
T_1 & \Diamond i \to i \\
\text{Density-in-itself} & \Diamond\neg i
\end{array}
$$
These properties are not definable in the basic modal language
(Corollary~\ref{c:t1-t6}).  $T_2$-separation, on the other hand,
remains undefinable even with nominals (Theorem~\ref{t:T_2
undefinable}).

Modal languages containing nominals are often called \emph{hybrid
  languages}.  In this section we investigate the topological
expressive power of two hybrid languages, namely $\mathcal{H}(@)$ and
$\mathcal{H}(E)$.  Formally, fix a countably infinite set of nominals
$\textsc{Nom}=\{i_1, i_2, \ldots\}$, disjoint from the set
$\textsc{Prop}$ of proposition letters. Then the formulas of
$\mathcal{H}(@)$ are given by the following recursive definition:
$$
\mathcal{H}(@) \qquad\qquad
  \phi  ~::=~  \top ~|~ p ~|~ i ~|~ \phi\wedge\phi ~|~ \neg\phi ~|~
  \Box\phi ~|~ @_i\phi
$$
where $p\in\textsc{Prop}$ and $i\in\textsc{Nom}$. $\mathcal{H}(E)$
further extends $\mathcal{H}(@)$ with the global modality, which was
described in the previous section. Thus, the formulas of
$\mathcal{H}(E)$ are given by the following recursive definition:
$$
  \qquad \mathcal{H}(E) \qquad\qquad
\phi ~::=~ \top ~|~ p ~|~ i ~|~ \phi\wedge\phi ~|~ \neg\phi ~|~
\Box\phi ~|~ @_i \phi ~|~ E\phi
$$
As in the previous section, we will use $A\phi$ as an abbreviation
for $\neg E \neg \phi$. We have introduced $@_i$ as a primitive
operator, but it will become clear after introducing the semantics
that $@_i$ can be defined in terms of the operator $E$.

\begin{defn}
  A {\em hybrid topological model} $\mathfrak{M}$ is a topological
  space $(X,\tau)$ and a valuation
  $\nu:\textsc{Prop}\cup\textsc{Nom}\to\wp(X)$ which sends
  propositional letters to subsets of $X$ and nominals to singleton
  sets of $X$.
\end{defn}
  The semantics for $\mathcal H(@)$ and $\mathcal H(E)$ is the same as for the basic
  modal language for the propositional letters, nominals, Boolean
  connectives, and the modality $\Box$. The semantics of $@$ and $E$
  is as follows:
\[\begin{array}{lll}
\mathfrak{M},w \models @_i \phi & \mathrm{iff} &
\mathfrak{M},v \models \phi \textrm{ for } \nu(i)=\{v\}\\
\mathfrak{M},w \models E\phi & \mathrm{iff} & \exists v \in X.
(\mathfrak{M},v \models \phi) \\
\end{array}\]
Validity and definability are defined as for the basic modal language,
but considering only valuations that assign singleton sets to the
nominals.

%\begin{defn}
%A relation $Z$ between hybrid models $\mathfrak M$ and $\mathfrak
%M'$ is called a \emph{hybrid topo-bisimulation} if it is a
%topo-bisimulation between the corresponding topological spaces and:
%
%\begin{description}
%\item[\textbf{Atom}] If $x Z x'$ then $x \in \nu(p)$ iff $x' \in \nu'(p)$ for
%  all $p \in \textsc{Prop} \cup \textsc{Nom}$
%\item[\textbf{Nom}] If $x \in \nu(i)$ and $x' \in \nu'(i)$ for
%  some $i \in  \textsc{Nom}$ then $x Z x'$
%\end{description}
%
%A topo-bisimulation $Z \subseteq X \times X'$ is called {\it total}
%if $\forall x {\in} X.\ \exists x' {\in} X'.(x Z x')$ and $\forall
%x' {\in} X'.\ \exists x {\in} X.(x Z x')$.
%
% We denote hybrid bisimulation between models $\mathfrak M$ and
% $\mathfrak M'$ as $\mathfrak M \bisim_{H(@)}\mathfrak M'$ and total
% hybrid bisimulation as $\mathfrak M \bisim_{H(E)} \mathfrak M'$.
%\end{defn}

\begin{prop}
  Taking topological sums or interior images in general does not preserve
  validity of $\mathcal{H}(@)$-formulas.
\end{prop}

\begin{proof} The one-point space $X=\{0\}$ with the trivial
  topology validates $@_i\Diamond j$, but this formula is
  not valid on $X \biguplus X$. Thus topological sums do
  not preserve validity for $\mathcal{H}(@)$.

  To see that $\mathcal{H}(@)$-validity is not preserved by
  interior maps, consider natural numbers with the topology induced by the ordering,
  i.e.\ the space $(\mathbb N, \tau)$ where $\tau = \{ [a, \infty) \mid
  a\in\N \}\cup\set{\emptyset}$. The formula $\varphi=@_i \Box(\Diamond i \to i)$ (which
  defines antisymmetry in the relational case) is easily seen to be
  valid in it. Then consider a topological space $X=\{0,1\}$ with the
  trivial topology $\tau'=\{\emptyset,X\}$ and a map $f$ that sends
  even numbers to $0$ and odd numbers to $1$. This is an interior map,
  however, $\varphi$ is not even satisfiable on $X$.
\end{proof}

On the other hand, the validity of $\mathcal{H}(@)$-formulas is
preserved under taking open subspaces.

\begin{lem} \label{lem:tufmi_preserv2} The validity of
  $\mathcal{H}(@)$ formulas is preserved under taking open subspaces.
\end{lem}

The proof is identical to that of Theorem~\ref{t:open
  subspaces}. Recall from Section~\ref{sec:open subspaces} that
connectedness is not preserved under taking open subspaces.  As a
corollary, we obtain that connectedness is not definable in
$\mathcal{H}(@)$.

Also, validity of $\mathcal{H}(E)$-formulas \emph{is} reflected by
Alexandroff extensions. We can in fact improve on this a bit, using
the notion of a \emph{topological ultrafilter morphic image}.

%\todo{``Topological ultrafilter morphic image'' sounds terribly complicated, can't we come up with a more catchy name?}

\begin{defn}
Let $X$ and $Y$ be topological spaces. $Y$ is called a {\em
topological ultrafilter morphic image} (or simply an \emph{u-morphic
image}) of $X$ if there is a surjective interior map $f: X \to
\Ae{Y}$ such that $|f^{-1}(\pi_y)| = 1$ for every principal
ultrafilter $\pi_y \in \Ae{Y}$ (one can say figuratively `$f$ is
injective on principal ultrafilters').
\end{defn}

Clearly, every space is a u-morphic image of its Alexandroff
extension.

\begin{lem} \label{lem:tufmi_preserv}
The validity of $\mathcal{H}(E)$ formulas is preserved under taking
\mbox{u-morphic} images.
\end{lem}

\proof Let $X$ and $Y$ be topological spaces and $f: X \to \Ae{Y}$
an interior map that is injective on principal ultrafilters. Suppose
further that $Y\not\models\phi$. We will show that $X
\not\models \phi$.

Since $Y\not\models\phi$, there is a valuation $\nu$ on $Y$ such that
$\nu(\phi)\neq Y$. Consider the valuation on $\Ae{Y}$ defined by
$\Ae{\nu}(p)=\{\mathfrak u\in Y^* \mid \nu(p)\in \mathfrak u\}$, where
$p$ can be a propositional letter or a nominal. It is not hard to see
that $\nu^*$ assigns to each nominal a singleton set consisting of a
principal ultrafilter. Next, we define the valuation $\nu'$ on $X$ by
$\nu'(p)= f^{-1}(\nu^*(p))$.  Since $f$ is injective on principal
ultrafilters, $\nu'$ again assigns singleton sets to the nominals.
Finally, a straightforward induction argument reveals that for all
$a\in X$ and $\psi\in \mathcal{H}(E)$,
$$(X,\nu'),a\models\psi\ \ \Leftrightarrow\ \ (Y^*,\nu^*),f(a)\models
\psi\ \ \Leftrightarrow\ \ \nu(\psi)\in f(a)$$
As $\nu(\neg\phi)\neq\emptyset$ there exists an ultrafilter
$\mathfrak u\in Y^*$ which contains $\nu(\neg\phi)$. Since $f$ is
onto, there exists $a\in X$ such that $f(a)=\mathfrak u$. It follows
that $(X,\nu'),a\models\neg\phi$ and therefore $X\not\models\phi$,
as required.
\endproof

The following two results characterize definability in
$\mathcal{H}(E)$ and $\mathcal{H}(@)$ in terms of closure under
taking u-morphic images. The proofs are inspired by relational
results presented in \cite{cate04}.

\begin{thm} \label{t:hybe_def} Let $\mathsf{K}$ be any
  $\mathcal{L}_t$-definable class of topological spaces. Then
  $\mathsf{K}$ is definable in $\mathcal{H}(E)$ iff $\mathsf{K}$ is closed under
  u-morphic images.
\end{thm}

\proof Lemma~\ref{lem:tufmi_preserv} constitutes the proof of the
left-to-right direction. We will prove the right-to-left direction.
Let $Log(\mathsf K)$ be the set of $\mathcal{H}(E)$-formulas valid on
$\mathsf K$. We will show that every space $X \models Log(\mathsf K)$
belongs to $\mathsf K$, and hence $Log(\mathsf K)$ \emph{defines}
$\mathsf{K}$.

Suppose $X\models Log(\mathsf K)$.  We introduce a propositional
letter $p_A$ for every subset $A \subseteq X$, as well as a nominal
$i_a$ for every $a \in X$. These propositional letters and nominals
are interpreted on $X$ by the natural valuation. Let $\Delta$ be the
following set of formulas, where $B$ and $C
$ range over all subsets
of $X$ and $a$ ranges over all points of $X$:
$$
\begin{array}{clc}
& A(i_a \leftrightarrow p_{\{a\}}) & \\
& A(p_{-B} \leftrightarrow \neg p_B) & \\
& A(p_{B \cap C} \leftrightarrow p_B \land p_C) & \\
& A(p_{\I B} \leftrightarrow \Box p_B) & \\
& A(p_{\C B} \leftrightarrow \Diamond p_B) & \\
\end{array}
$$

As in the proof of Theorem~\ref{t:GbTh global}, we can find an
$\mathcal{L}_t$-saturated (hybrid) topological model, based on a
space $Y\in\mathsf K$, that makes $\Delta$ globally true.  Note that
conditions (1)--(5) from the proof of Theorem~\ref{t:up_ae} hold for
$Y$ (the truth of $A(i_a \leftrightarrow p_{\{a\}})$ ensures that
the predicates $P_{\{a\}}$ have non-empty interpretation). It
follows, by the same argument as in the proof of
Theorem~\ref{t:up_ae}, that the map $f: Y \to \Ae{X}$ defined by
$$
f(a)=\{A \subseteq X \mid \textstyle Y \models
P_A(a) \}
$$
is a surjective interior map.  We will now show that $f$ is
injective on principal ultrafilters. Suppose there exist $w, v \in
Y$ and $f(w)=f(v)=\pi_a$ where $a \in X$ and $\pi_a$ is the
principal ultrafilter containing $\{a\}$. By definition of $f$ we
get $Y \models P_{\{a\}}(w)\wedge P_{\{a\}}(v)$.  By global truth of
$\Delta$ we obtain $Y, w \models i_a$ and $Y, v \models i_a$, hence
$w=v$.

It follows that $X$ is an u-morphic image of $Y$. As $\mathsf K$ is
closed under u-morphic images, we conclude that $X\in\mathsf K$ as
required. \endproof

\begin{thm} \label{t:hybat_def}
Let $\mathsf{K}$ be any $\mathcal{L}_t$-definable class of
  topological spaces. Then $\mathsf{K}$ is definable in $\mathcal H(@)$ iff it is closed under
topological ultrafilter morphic images and under taking open
subspaces.
\end{thm}

\proof The `only if' part is taken care of by
Lemmata~\ref{lem:tufmi_preserv2} and \ref{lem:tufmi_preserv}.  The
proof of the `if' part proceeds as in Theorem~\ref{t:GbTh basic Lt},
with some modifications.

The first difference is that the set of formulas $\Delta$ is augmented
with formulas of the form $@_{i_a} p_{A}$, for all points $a$ that
belong to a non-empty set $A\subseteq X$.

A compactness argument similar to the one used in the proof of
Theorem~\ref{t:GbTh basic Lt} shows that  $\{@_{i_a}\Box\phi \mid
\phi\in\Delta, a\in X\}$ is true in some $\mathcal{L}_t$-saturated
topological model $\mathfrak{N}=(Y,\mu)$ with $Y\in\mathsf{K}$. For
each $b\in Y$ named by a nominal, choose an open neighborhood $O_b$
as described in Definition~\ref{def:Lt-saturation}. Let $O$ be the
union of all these open neighborhoods. Note that by closure under
open subspaces we obtain $O\in\mathsf K$. It is not hard to see that
the submodel $\mathfrak{K}$ of $\mathfrak{N}$ based on the open
subspace $O$ globally satisfies $\Delta$, and hence satisfies the
conditions (1)--(5) described in the proof of Theorem~\ref{t:up_ae}.

Thus there exists an interior map $f$ from $O$ onto $\Ae{X}$. That
$f$ is injective on principal ultrafilters can be proved as in
Theorem~\ref{t:hybe_def}. Thus $X$ is an u-morphic image of
$O\in\mathsf K$. Since $\mathsf K$ is closed under u-morphic images,
we obtain $X\in\mathsf K$ as required.
\endproof

As an application, we will show that $\mathcal{H}(@)$ and
$\mathcal{H}(E)$ are not expressive enough to be able to define the
$T_2$ separation property. Recall the definition of irresolvability
(Definition~\ref{def:irresolvability}). We call a space $X$
\emph{$\alpha$-resolvable} for a cardinal number $\alpha$ if $X$
contains $\alpha$-many pairwise disjoint dense subsets. In
\cite{feng}, an $2^{2^{\aleph_0}}$-resolvable $T_2$-space was
constructed. We use this space to prove that

\begin{thm}\label{t:T_2 undefinable}
  The class of $T_2$ topological spaces in not definable in
  $\mathcal{H}(@)$ and $\mathcal{H}(E)$.
\end{thm}

\proof We employ an argument similar to, but more complicated than,
the one used in Corollary~\ref{c:t1-t6}. Our strategy is as follows:
we construct spaces $X$ and $Y$ such that: $Y$ is a $T_2$ space, $X$
is an u-morphic image of $Y$, and $X$ is not a $T_2$ space.  Then we
apply Theorem~\ref{t:hybe_def}.

%Since $T_1$ can be expressed in $H(@)$ and the validity of hybrid
%formulas is preserved under u-morphic images, a space $X$ should
%satisfy $T_1$ but should not satisfy $T_2$. We proceed to describe
%such a space.

Take $X=(\mathbb{N},\tau)$ where $\tau$ is the co-finite topology.
That is
$$\tau=\set{\emptyset}\cup \set{A\subseteq\N ~\mid~ \N{\setminus}A \textrm{ is finite}}$$
Then $X$ is $T_1$ since every singleton is closed, but not $T_2$ as
any two non-empty opens necessarily meet.
Denote by $\mathfrak F$ the set of all the free ultrafilters over
$\N$. Then the following holds:
\begin{quote}%\label{c:describe X*}
  \emph{Claim 1:} The topology $\tau^*$ of the Alexandroff extension
  $\Ae{X}$ of $X$ is described as follows:
$$O\in\tau^*\ \ \ \ \Leftrightarrow\ \ \ \ \ \mathfrak F\subseteq
O$$

\medskip\emph{Proof:}
Suppose $O\in\tau^*$. If $O=X^*$ the claim follows.
Otherwise $O$ contains a basic open set $\mathcal G^*$ which
consists of all the ultrafilters extending a proper open filter
$\mathcal G$. Note that if $A\in\mathcal G$ is not cofinite, then
$\I A=\emptyset\notin\mathcal G$. Therefore, $\mathcal G$ consists
of cofinite sets only. Since each free ultrafilter contains
\emph{all} cofinite sets, we obtain $\mathfrak F\subseteq \mathcal
G^*\subseteq O$.

Now for the other direction. Suppose $\mathfrak F\subseteq O$. First
note that $\mathfrak F$, being the extension of the open filter of
all cofinite subsets of $X$, is a basic open in $\tau^*$. Further,
if $x\in X$, then the open filter $O_x=\set{A\mid x\in A,\ A\textrm
{ cofinite}}$ is such that $O_x^*=\pi_x\cup F_x$ where $\pi_x$
denotes the principal filter of $x$ and $F_x\subseteq \mathfrak F$.
It follows that
$$O=\mathfrak F \cup \bigcup\limits_{\pi_x\in O}(\pi_x\cup F_x)$$
Since each $\pi_x\cup F_x=O_x^*\in\tau^*$ we obtain that
$O\in\tau^*$. The claim is proved.
\end{quote}
Next we will construct the space $Y$. Let $Z=(Z,\tau_1)$ be a
$2^{2^{\aleph_0}}$-resolvable topological space which satisfies
$T_2$ (according to \cite{feng} such a space exists). We will denote
$2^{2^{\aleph_0}}$ many dense disjoint subsets of $Z$ by $Z_\iota$
where $\iota \in \mathfrak F$. Here $\mathfrak F$ is again the set
of all free ultrafilters over $\mathbb{N}$. Since the cardinality of
$\mathfrak F$ is known to be $2^{2^{\aleph_0}}$ \cite[Corollary
3.6.12]{engelking}, such indexing is possible. Let $\bar{Z} =
Z-\bigcup\limits_{\iota
  \in \mathfrak F} Z_\iota$.  Thus
$$Z = \bar{Z}\cup \bigcup\limits_{\iota \in \mathfrak F} Z_\iota$$

Put $Y=(\mathbb{N} \cup Z, \tau')$ where $\tau'$ is as follows:
$$\tau'=\set{\emptyset}\cup \set{O\subseteq Y \mid O\cap Z\in\tau_1, O\cap Z\neq\emptyset}$$

In words---the topology of $Z$ as a subspace of $Y$ is $\tau_1$ and
the neighborhoods of the points from $\N$ are the sets of the form
$\{x\} \cup O$ where $x \in \mathbb{N}$ and $\emptyset\neq
O\in\tau_1$.

\begin{quote}%\label{c:Y is T_2}
\emph{Claim 2:} $Y$ is a $T_2$ space.

\medskip\emph{Proof:} Indeed, any two points that belong to $Z$ can be
separated by two opens from $\tau_1$, since $(Z,\tau_1)$ is a $T_2$
space. Any two points $x, y \in \mathbb{N}$ can be separated by open
sets of the form $\{x\} \cup O_x$ and $\{y\} \cup O_y$ where $O_x,
O_y\in\tau_1$ are non-empty open sets from $Z$ such that $O_x \cap O_y
= \emptyset$. Finally, two points $x,y$ such that $x \in \mathbb{N}$
and $y \in Z$ can be separated by the sets $\{x\} \cup O_x$ and $O_y$
where again $O_x$ and $O_y$ are disjoint non-empty open subsets of
$Z$.
\end{quote}
Now we construct the mapping $f: Y\to X^*$.
Pick any $\zeta\in\mathfrak F$ and define $f:
\mathbb{N} \cup Z \to X$ as follows:
$$
f(x) = \left\{\begin{array}{ll}
\pi_x & \textrm{ if } x \in \mathbb{N}\\
\iota & \textrm{ if } x \in Z_\iota\\
\zeta & \textrm{ if } x \in \bar{Z}
\end{array}\right.
$$

\begin{quote}
\emph{Claim 3:} The map $f$ is a surjective interior map.

\medskip\emph{Proof:} That $f$ is surjective follows from the
construction.

Let us show that $f$ is continuous. Take $O\in\tau^*$. By
Claim~1 we have $\mathfrak F\subseteq O$. It
follows from the definition of $f$ that $f^{-1}O$ is of the form $Z
\cup A$ where $A \subseteq \mathbb{N}$. From the definition of
$\tau'$ we obtain $f^{-1}(O)\in\tau'$.

To show that $f$ is an open map, take an arbitrary open set $O \in
\tau'$. It follows from the definition of $\tau'$ that $O\cap
Z\in\tau_1$ and $O\cap Z\neq \emptyset$. Then, as each $Z_\iota$ is
dense in $Z$, it follows that $O\cap Z_\iota\neq\emptyset$ for all
$\iota \in \mathfrak F$. Hence, $f(O)$ contains $\mathfrak F$ and is
open in $X^*$ according to Claim~1.
\end{quote}
Note that $f$ is injective on principal ultrafilters, by
construction. Therefore $X$ is an u-morphic image of $Y$. Since $Y$
is $T_2$ and $X$ is not, it follows that the class of $T_2$ spaces
is not closed under u-morphic images. Recall that the class of $T_2$
spaces is $\mathcal{L}_t$-definable. It follows by
Theorem~\ref{t:hybe_def} and Theorem~\ref{t:hybat_def} that the
class of $T_2$ spaces is not definable in $\mathcal{H}(E)$ and
$\mathcal{H}(@)$.
\endproof

%}}}

%{{{ Language extensions (difference mod-ty, binder)
\subsection{The difference modality}

In this section, we consider $\mathcal{M}(D)$, the extension of the
basic modal language with the \emph{difference modality} $D$.  Recall
that the global modality allows us to express that a formula holds
\emph{somewhere}. The \emph{difference modality} $D$ allows us to
express that a formula holds \emph{somewhere else}. For example,
$p\land\neg Dp$ expresses that $p$ is true at the current point and
nowhere else. Formally,
$$
\mathfrak M, w \models D \varphi\ \ \ \ \mbox{ iff }\ \ \ \ \exists
v \neq w. (\mathfrak M, v \models \varphi)
$$
The global modality is definable in terms of the difference
modality: $E\phi$ is equivalent to $\phi\lor D\phi$.  It follows
that $\mathcal{M}(D)$ is at least as expressive as $\mathcal{M}(E)$.
Furthermore, one can express in $\mathcal{M}(D)$ that a
propositional letter $p$ is true at a unique point (i.e., behaves as
a nominal): this is expressed by the formula $E(p\land\neg Dp)$.
Combining these two observations, it is not hard to show that every
class of topological spaces definable in $\mathcal{H}(E)$ is also
definable in $\mathcal{M}(D)$. The opposite also holds
\cite{gargov93,tadeusz}:

\begin{thm}
  $\mathcal{M}(D)$ can define exactly the same classes of topological
  spaces as $\mathcal{H}(E)$.
\end{thm}

\begin{cor}\label{c:GbTh_D}
  An $\mathcal L_t$-definable class of topological spaces is definable
  in $\mathcal{M}(D)$ iff it is closed
  under u-morphic images.
\end{cor}
Recall that the separation axioms $T_0$ and $T_1$, as well as
\emph{density-in-itself}, are definable in the language
$\mathcal{H}(E)$.  They are definable in $\mathcal{M}(D)$ as
follows, where $U\phi$ is short for $\phi\land\neg D\phi$:
$$
\begin{array}{ll}
T_0:       &  U p\wedge DU q\to\Box\neg q\vee D( q\wedge\Box\neg p) \\
T_1:       & U p\to A( p\leftrightarrow\Diamond p) \\
\text{Density-in-itself}: \quad & p \to \Diamond Dp
\end{array}
$$
For more on topological semantics of $\mathcal{M}(D)$ we refer to a
recent study \cite{kudinov}.

\subsection{The $\downarrow$-binder}

The last extension we will consider is the one with explicit point
variables, and with the $\downarrow$-binder. The point variables are
similar to nominals, but their interpretation is not fixed in the
model. Instead, they can be bound to the current point of evaluation
using the $\downarrow$-binder.  For instance, $\downarrow\! x.\Box x$
expresses that the current point is an isolated point.

%\todo{Can we find a natural example of a class of spaces that becomes
%  definable with $\downarrow$?}

$\mathcal{H}(@,\downarrow)$ and $\mathcal{H}(E,\downarrow)$ are the
extensions of $\mathcal{H}(@)$ and $\mathcal{H}(E)$, respectively,
with state variables and the $\downarrow$-binder.  Formally, let
$\textsc{Var}=\{x_1, x_2, \ldots\}$ be a countably infinite set of
point variables, disjoint from \textsc{Prop} and \textsc{Nom}.  The
formulas of $\mathcal{H}(@,\downarrow)$ and
$\mathcal{H}(E,\downarrow)$ are given by the following recursive
definitions (where $p\in\textsc{Prop}$, $i\in\textsc{Nom}$, and
$x\in\textsc{Var}$):
\[\begin{array}{ll}
\mathcal{H}(@,\downarrow) \quad &
  \phi ~::=~ p ~\mid~ i ~\mid~ x ~\mid~ \neg\phi ~\mid~ \phi\land\psi
  ~\mid~ \Box\phi ~\mid~ @_i\phi ~\mid~ \downarrow\! x.\phi \\[2mm]
\mathcal{H}(E,\downarrow) \quad &
  \phi ~::=~ p ~\mid~ i ~\mid~ x ~\mid~ \neg\phi ~\mid~ \phi\land\psi
  ~\mid~ \Box\phi ~\mid~ @_i\phi ~\mid~ E\phi ~\mid~ \downarrow\! x.\phi
\end{array}\]
These formulas are interpreted, as usual, in topological models.
However, the interpretation is now given relative to an
\emph{assignment} $g$ of points to point variables (just as in
$\mathcal{L}_t$). The semantics of the state variables and
$\downarrow$-binder is as follows:
\[\begin{array}{lll}
  \mathfrak{M},w,g\models x &\text{iff}& g(x)=w \\
  \mathfrak{M},w,g\models \downarrow\! x.\phi &\text{iff}&
    \mathfrak{M},w,g^{[x\mapsto w]}\models\phi
  \end{array}\]
  where $g^{[x\mapsto w]}$ is the assignment that sends $x$ to $w$ and
  that agrees with $g$ on all other variables.
  We will restrict attention to \emph{sentences}, i.e., formulas in
  which all occurences of point variables are bound. The
  interpretation of these formulas is independent of the assignment.

  It turns out that $\mathcal H(E,\downarrow)$ is essentially a notational
  variant for a known fragment of $\mathcal{L}_t$, called
  $\mathcal{L}_I$. This is the fragment of $\mathcal{L}_t$
  where quantification over opens is only allowed in the form,
  for $U$ not occurring in $\alpha$:
\[\begin{array}{ll}
  \exists U.(x\varepsilon U \land \forall y.(y\varepsilon U\to
  \alpha)),
  & \text{ abbreviated as $[I_y\alpha](x)$, and, dually,} \\
  \forall U.(x\varepsilon U \to \exists y.(y\varepsilon U\land
  \alpha)),
  & \text{ abbreviated as $[C_y\alpha](x)$}.
\end{array}\]

Comparing the above with the Definition~\ref{d:ST} reveals that the
formulas of the basic modal language translate inside
$\mathcal{L}_I$ by the standard translation. So $\mathcal{ML}$ can
be thought of as a fragment of $\mathcal{L}_I$. Apparently, adding
nominals, $\downarrow$ and $E$ to the language is just enough to get
the whole of $\mathcal{L}_I$.

\begin{thm} $\mathcal H(E,\downarrow)$ has the same expressive power as
  $\mathcal{L}_I$.
\end{thm}

\begin{proof}
  The standard translation from modal logic to $\mathcal{L}_t$ can be
  naturally extended to $\mathcal{H}(E)$, treating nominals  as first-order constants. The extra clauses are then

  \[\begin{array}{lll}
     ST_x(t) &=& x=t \quad \text{for $t\in\textsc{Nom}\cup\textsc{Var}$} \\
     ST_x(@_i\varphi) &=& \exists x.(x=c_i \land ST_x(\varphi)) \\
     ST_x(E\varphi) &=& \exists x.ST_x(\varphi) \\
     ST_x(\downarrow y.\varphi) &=& \exists y.(y=x\land ST_x(\varphi))
  \end{array}\]
  It is easily seen that this extended translation maps
  $\mathcal{H}(E,\downarrow)$-sentences to $\mathcal{L}_I$-formulas in
  one free variable.
  Conversely, the translation $HT_x$ below maps
  $\mathcal{L}_I$-formulas $\alpha(x)$ to
  $\mathcal{H}(E,\downarrow)$-sentences:

  \[\begin{array}{lll}
     HT(s=t) &=&  @_s t \\
     HT(Pt)  &=&  @_t p \\
     HT(\neg\alpha) &=& \neg HT(\alpha) \\
     HT(\alpha\land\beta) &=& HT(\alpha)\land HT(\beta) \\
     HT(\exists x.\alpha) &=& E\downarrow\! x.HT(\alpha) \\
     HT([I_y\alpha](t)) &=& @_t \Box \downarrow y.HT(\alpha) \\
     HT([C_y\alpha](t)) &=& @_t \Diamond \downarrow y.HT(\alpha) \\[2mm]
     HT_x(\alpha(x)) &=& \downarrow\! x.HT(\alpha)
  \end{array}\]
  It is not hard to see that both translations preserve truth,
  in the sense of Theorem~\ref{t:stand-transl}.
\end{proof}
This connection allows us to transfer a number of known results.
For instance, $\mathcal{L}_I$ has a nice axiomatization, it is know
to have interpolation, and the $\mathcal{L}_I$-theory of the class
of $T_1$-spaces is decidable (see \cite{makowsky}). Hence, these
results transfer to $\mathcal H(E,\downarrow)$.  It is also known
that $\mathcal{L}_I$ is strictly less expressive than
$\mathcal{L}_t$. In particular, there is no $\mathcal{L}_I$-sentence
that holds precisely on those topological models that are based on a
$T_2$-space. Hence, the same
holds for $\mathcal{H}(E,\downarrow)$.%
\footnote{In fact, Makowsky and Ziegler \cite{makowsky} showed that,
in the absence of proposition letters and nominals, every two
dense-in-itself $T_1$-spaces have the same $\mathcal{L}_t$-theory.}
Note that this does not imply \emph{undefinability} of $T_2$ in
$\mathcal{H}(E,\downarrow)$. Nevertheless, we conjecture that $T_2$
is not definable in $\mathcal{H}(E,\downarrow)$.

The precise expressive power of $\mathcal{L}_I$ on topological models
can be characterized in terms of \emph{potential homeomorphisms}.

\begin{defn}
  A \emph{potential homeomorphism} between topological models
  $\mathfrak{M}=(M,\tau,\nu)$ and $\mathfrak{N}=(N,\sigma,\mu)$ is a
  family $F$ of  partial bijections
  $f:M\to N$ satisfying the following conditions for each $f\in F$:
\begin{enumerate}
\item $f$ preserves truth of proposition letters and nominals
  (in both directions).
\item - For each $m\in M$ there is a $g\in F$ extending
  $f$, such that $m\in dom(g)$.

  - For each $n\in N$, there is a $g\in F$ extending $f$,
  such that $n\in rng(g)$.
\item - For each $(m,n)\in f$ and open neighborhood $U\ni m$,
  there is an open neighborhood $V\ni n$ such that
  for all $n'\in V$ there is a $g\in F$ extending $f$ and an $m'\in U$
  such that $(m',n')\in g$.

  - Likewise in the opposite direction.
\end{enumerate}
\end{defn}

The following characterization follows from results in
\cite{makowsky}.

\begin{thm} An $\mathcal{L}_t$-formula $\phi(x_1, \ldots, x_n)$ is
  equivalent to an $\mathcal{L}_I$-formula in the same free variables
  iff it is invariant for potential homeomorphisms.
\end{thm}

Substituting $\mathcal{H}(E,\downarrow)$ for $\mathcal{L}_I$, this
gives us a Van Benthem-style characterization of
$\mathcal{H}(E,\downarrow)$ as a fragment of $\mathcal{L}_t$.  We
leave it as an open problem to find a similar characterization of
the expressive power of $\mathcal H(@,\downarrow)$. We also leave it
as an open problem to characterize the classes of topological spaces
definable in these languages.

Note that the union of the graphs of the partial bijections that
constitute a potential homeomorphism gives rise to a total
topo-bisimulation between the models in question. Thus a formula
that is invariant for topo-bisimulations is also invariant for
potential homeomorphisms. This is a semantical side of the fact that
the basic modal language is a fragment of $\mathcal{L}_I$. In fact,
we could have taken the language $\mathcal{L}_I$ as our first-order
correspondence language from the very beginning. A feeling that
$\mathcal{L}_I$ might be the `right' candidate for the topological
correspondence language might be strengthened by the fact that in
its relational interpretation (i.e., on Kripke structures),
$\mathcal H(E,\downarrow)$ has the full expressive power of the
first-order correspondence language. We stand, however, by our
choice of $\mathcal{L}_t$ since: (a) it provides stronger
definability results (there are more $\mathcal{L}_t$-definable
classes than $\mathcal{L}_I$-definable ones); (b) $\mathcal{L}_t$ is
closer to both the usual first-order signature and the usual
set-theoretic language used to formalize concepts in general
topology.

%}}}

%{{{ Applications (separation axioms etc.)

%}}}

%{{{ Further work

\section{Discussion}\label{s:discussion}

%\todo{Expand and polish this section}

We have studied the expressive power of various (extended) modal
languages interpreted on topological spaces.
Tables~\ref{tab:overview-GbTh} and~\ref{tab:examples} summarize and
illustrate our main findings, concerning definability of classes of
spaces. We also obtained a Van Benthem-style characterization of the
basic modal language in terms of topo-bisimulations, thereby solving
an open problem from \cite{blackburnvanbenthem}. 

Some of the key innovative elements in our story are (\emph{i})
identifying the appropriate topological analogues of familiar
operations on Kripke frames such as taking bounded morphic images, or,
ultrafilter extensions (\emph{ii}) identifying $\mathcal{L}_t$ as
being the appropriate correspondence language on topological models
(indeed, our result confirm once again that, as has been claimed
before, $\mathcal L_t$ functions as the same sort of
``landmark'' in the landscape of topological languages as first-order
logic is in the landscape of classical logics), and (\emph{iii})
formulating the right notion of saturation for $\mathcal{L}_t$ (which
many of our technical proofs depend on).

Our results on the hybrid language $\mathcal{H}(E,\downarrow)$ are
remarkable. For example, they show that, while
$\mathcal{H}(E,\downarrow)$ is expressively equivalent to the
first-order correspondence language on relational structures, it is
strictly less expressive than $\mathcal{L}_t$ on topological
models. This seems one more instance of the more sensitive power of
topological modeling.

Given that Alexandroffness is definable in $\mathcal{L}_t$, many of
our results can be seen as generalizing known results for modal languages
on (transitive reflexive) relational structures, and it is quite well
possible that results on the topological semantics will yield new
consequences for the relational semantics.

\begin{table}
\caption{Definability in extended modal languages}
\label{tab:overview-GbTh}{\small
\begin{center}
\begin{tabular}{l|p{.65\linewidth}|l}
  & characterization of definability for $\mathcal{L}_t$-definable classes & \\[2mm]
\hline \\[-2mm]
  $\mathcal{ML}$ & \emph{closed under topological sums, open subspaces and interior images, reflecting Alexandroff extensions} & Theorem~\ref{t:GbTh basic Lt} \\[6mm]
  $\mathcal{M}(E)$ & \emph{closed under interior images, reflecting Alexandroff extensions} & Theorem~\ref{t:GbTh global}  \\[6mm]
  $\mathcal{H}(@)$ & \emph{closed under open subspaces and u-morphic images}  & Theorem~\ref{t:hybat_def} \\[6mm]
  $\mathcal{H}(E)$ & \emph{closed under u-morphic images} & Theorem~\ref{t:hybe_def}\\ [6mm]
  $\mathcal{M}(D)$ & \emph{closed under u-morphic images} &
  Corollary~\ref{c:GbTh_D}
\end{tabular}
\end{center}}
\end{table}

\begin{sidewaystable}\small
\begin{center}
\caption{Properties of topological spaces definable in various
languages} \label{tab:examples}

\vspace{9mm}

\begin{tabular}{c|c|cccccc}
    & $\mathcal{L}_t$ & $\mathcal{M}$ & $\mathcal{M}(E)$ & $\mathcal{H}(@)$ & $\mathcal{H}(E)$ (or $\mathcal{M}(D)$)  \\
\hline &&&&& \\[5mm]
$T_0$ & \begin{tabular}{c}$\forall xy.(x\neq y\to$\\$\exists U_x.(y\nin U)\vee \exists V_y.(x\nin V))$\end{tabular}  & no & no & \begin{tabular}{c}$@_i\Diamond j \land @_j\Diamond i$\\$\to @_ij$\end{tabular} & idem \\[5mm]
$T_1$ & $\forall xy.(x\neq y\to \exists U_x.(y\nin U))$ & no & no & $\Diamond i \to i$ & idem \\[5mm]
$T_2$ & \begin{tabular}{c}$\forall xy.(x\neq y\to$\\$\exists U_x.\exists V_y.\forall z.(z\nin U\vee z\nin V))$\end{tabular} & no & no & no & no \\[5mm]
\begin{tabular}{c}Density-in-\\itself\end{tabular} & $\forall x\forall U_x(\exists y\neq x(y\in U))$ & no & no & $\Diamond\neg i$ & idem \\[5mm]
Connectedness & no & no & \begin{tabular}{c}$\mathsf{A}(\Box p \lor \Box\neg p)$\\$\to \mathsf{A}p\lor\mathsf{A}\neg p$\end{tabular} & no & as in $\mathcal{M}(E)$ \\[5mm]
\begin{tabular}{c}Hereditary\\irresolvability\end{tabular} & no & \begin{tabular}{c}$\Box(\Box(p\to\Box p)\to p)$\\$\to \Box p$\end{tabular} & idem & idem & idem  \\
\end{tabular}
\end{center}
\end{sidewaystable}

We finish by mentioning interesting directions for future research.

\begin{itemize}
\item {\bf Correspondence theory for alternative semantics.} There
are at least two other semantic paradigms where the approach taken
in this paper might prove useful. We discuss them briefly.

{\bf Diamond as derived set operator}. For any
      subset $S$ of a topological space, the derived set ${\mathbf d} S$ is the
      set of limit points of $S$, i.e., all points $x$ of which each
      open neighborhood contains an element of $S$ distinct from $x$
      itself. The closure operator can be defined in terms of the
      derived set operator: $\mathbb{C}S = S\cup {\mathbf d}S$. The converse
      does not hold, as ${\mathbf d}$ is strictly more expressive than
      $\mathbb{C}$ \cite{shehtman, submax}. Indeed, if we interpret the $\blacklozenge$ as the
      derived set operator, then the modal formula $\blacklozenge\top$
      defines density-in-itself. With the derived set operator we can
      also partially mimic nominals: $p\land\neg\blacklozenge p$
      expresses that, within small enough neighborhoods, $p$ acts as a
      nominal for the current point. Conversely, with nominals we can
      partially mimic the ${\mathbf d}$-operator: $@_i(\blacklozenge\phi
      \leftrightarrow \Diamond(\phi\land\neg i))$ is valid. The precise
      connection between ${\mathbf d}$ and nominals remains to be
      investigated. The standard translation should be modified in the following way to account for the new semantics:
       $$ST_x(\blacklozenge\phi)~:=~\forall U.(x\varepsilon U\to\exists y.(y\varepsilon U\wedge x\neq y\wedge
       ST_y(\phi)))$$
This is still a $\mathcal L_t$-formula. Whether or not the
expressive power of the corresponding fragment of $\mathcal L_t$ can
be characterized in a way we have presented here remains to be seen.
The interested reader is referred to \cite{Es81, shehtman, submax}
for more details on this topological semantics.

    {\bf Neighborhood semantics.} This is a generalization of the
    topological semantics for modal logic that allows to tackle
    non-normal modal logics. The corresponding structures are
    \emph{neighborhood frames} $(W, n)$ where $W$ is a non-empty
    set and $n\subseteq W{\times}\wp(W)$ is a binary relation between points of $W$ and
    subsets of $W$. The correspondence theory for \emph{Monotonic} modal logic
    has been explored in \cite{hhansen} and analogues of
    the Goldblatt-Thomason theorem and the Van Benthem theorem have
    been proved for the neighborhood semantics. Quite general definability
    results are also put forth in \cite{kurz} using the
    co-algebraic approach. We believe that a modification of the
    approach presented in this paper will further strengthen the
    investigation of the precise expressive power of non-normal
    modal logics over neighborhood frames. We outline one
    possible route in this direction that is similar to, but more general than, the one pursued in
    \cite{hhansen}.

    Extend the language $\mathcal L^2$ by another
    intersorted binary relation symbol $\eta$. To ensure that $\varepsilon$ behaves like
    the membership relation we postulate
$$\forall U,V.(U=V\leftrightarrow
        \forall x.(x\varepsilon U\leftrightarrow x\varepsilon V))$$
    The models for this new language $\mathcal L^2_\eta$ are of the form
    $(X,\sigma,\nu)$ with $X$ a set and $\sigma\subseteq\wp X$. The
    relation $\varepsilon$ is interpreted as set-theoretic
    membership, while the interpretation of $\eta$ defines a
    relation between elements of $X$ and elements of $\sigma$ so
    that $(X,\eta^X)$ becomes a neighborhood frame. Conversely, each
    neighborhood frame $(W,n)$ gives rise to a structure $(W,\set{A\subseteq W\mid\exists w\in W.(wnA)})$.

    The standard translation can also be modified to suit the new semantics:
\[
    \begin{array}{rcl}
    ST_x(\Box\phi)&~:=~&\exists U.(x\eta U\wedge \forall y.(y\varepsilon U\to
    ST_y(\phi)))
%    ST_x(\Diamond\phi)&~:=~&\neg\exists U.(x\eta U\wedge \forall y.(y\nin U\to
%    ST_y(\phi)))
    \end{array}
\]
    Note that the whole story now becomes simpler than in the
    case of topological semantics, since there is no restriction on
    $\eta$. Thus the full apparatus of the model theory for first-order
    logic is at hand when considering the expressivity and
    characterization of modal logic over neighborhood frames, in
    terms of $\mathcal L^2_\eta$. At least this is the case for the
    modal logic $\mathbf E$ determined by the class of \emph{all}
    neighborhood frames. The situation might change if some other non-normal modal logic
    is taken as a base. We briefly discuss two representative examples.

    Consider the modal logic determined by the class of neighborhood frames that are closed under
    intersection. The closure under intersection is easily seen
    to be $\mathcal L^2_\eta$-definable by the formula:
    $$\forall x.\forall U,V.\left[x\eta U\wedge y\eta V\to \exists W.(y\eta W\wedge\forall z.(z\varepsilon U\wedge z\varepsilon V\leftrightarrow z\varepsilon W) )\right]$$
    Consequently, we expect the situation in this and similar, $\mathcal L^2_\eta$-definable cases
    to be rather straightforward.

    However, consider the modal logic $\mathbf M$,
    determined by the class of all \emph{monotone} neighborhood
    frames. Recall that $(W, n)$ is monotone, if $wnA$ and
    $A\subseteq B$ imply $wn B$. This condition is not expressible
    in $\mathcal L^2_\eta$, so in this case part of the story we
    witnessed in this paper might re-appear. That is to say, one
    needs to find a well-behaved fragment of $\mathcal L^2_\eta$
    that is invariant for monotone frames. One possibility is to
    restrict the quantification over open variables by admitting only formulas of the kind
    $\exists U.(x\varepsilon U\to \phi)$ with $U$ occuring positively in
    $\phi$. We leave it to further research to decide whether fully developing
    this approach is worthwhile in this and other interesting cases.
  \item {\bf Further extensions of the language}. One could consider
    other extensions of the modal language, e.g., with propositional
    quantifiers \cite{kremer} or fixed point operators
    \cite{vanbenthemsarenac}. It seems worthwhile to consider the
    extension of the signature of $\mathcal L_t$ with \emph{function}
    symbols (to model continuous transformations of spaces) or a
    \emph{binary relation} symbol (to model the time flow) and
    consider the applications to the domain of Dynamic Topological
    Logics of \cite{artemov97modal, Kremer-mints05} or other
    structures for modal spatio-temporal logics.

\item {\bf Axiomatizations}. In this paper, we have investigated
  \emph{expressive power} of extended modal languages interpreted on
  topological spaces. However, in order for these logics to be of
  practical use, their \emph{proof theory} will have to be studied as
  well. In the case of the basic modal language, topological
  completeness has already been studied for a long time \cite{mckinsey44},
  but for more expressive modal languages, this is a new area of
  research. Some first results for $\mathcal{M}(D)$ and
  $\mathcal{H}(E)$ and can be found in \cite{kudinov,tencatelitak}.

%  \todo{Is the above correct?
%        Is there anything more to say?

%        Just a loose thought: do ultrafilter
%        morphisms preserve validity of formulas with $d$? Probably not... }

%\item Another direction is to consider modal languages for more
 %     concrete spaces, e.g., metric spaces.

%  \todo{I don't know much about this, but i guess there is quite a
%      lot of recent research on this, perhaps we can give some
%      pointers?}
\end{itemize}

% \paragraph{Kripke semantics revisited}
%
% Recall that Alexandroffness is definable in $\mathcal{L}_t$.  If we
% resrict attention to Alexandroff spaces, and if we consider their
% relational representation, we see the following: the operations of
% topological sum, open subspace, interior image, and Alexandroff
% extension correspond to the relational operations of disjoint union,
% generated subframe, bounded morphic image, and ultrafilter extension.
% Hence, our topological version of the Goldblatt-Thomason theorem
% really \emph{generalizes} the relational result, since it implies the
% latter, at least on reflexive transitive frames.

%}}}

\bibliographystyle{plain}
\bibliography{tengasus}

\pagebreak

\appendix

\section{$\mathcal{L}^2$ over topological models}
\label{app:L2}

In this section, we prove Theorem~\ref{thm:L2-is-bad}, here
stated once again for reference:

\begin{thm} $\mathcal{L}^2$ interpreted on topological models
  lacks Compactness, L\"owenheim-Skolem and Interpolation, and
  is $\Pi^1_1$-hard for validity.
\end{thm}

This was already known for the more general case where
$\mathcal{L}^2$-formulas can contain $k$-ary relation symbols with
$k\geq 2$. The topological models we work with in this paper contain
only unary predicates, but we will show that the bad properties of
$\mathcal{L}^2$ already occur in this more restricted setting.

\begin{proof}
  These facts can all be derived from the observation that
 $\mathcal{L}^2$ can define $(\mathbb{N},\leq)$ up to
  isomorphism.

\begin{itemize}
\item Definability of $(\mathbb{N},\leq)$.

  Let $x\leq y$ stand for the $\mathcal{L}^2$-formula $\forall
  U.(x\varepsilon U \to y\varepsilon U)$, which defines the well known
  \emph{specialisation order} ($x\leq y$ iff $x\in\mathbb{C}\{y\}$).
  For each topological space $(X,\tau)$, $\leq$ defines a quasi-order
  on $X$.  Conversely, every quasi-order on a set $X$ is the
  specialisation order of some topology on $X$ (in fact, of an
  Alexandroff topology on $X$).

  A special feature of $\leq$ is that every open set $U$ is an
  \emph{up-set} with respect to $\leq$ (i.e., whenever $x\varepsilon
  U$ and $x\leq y$ then also $y\varepsilon U$).  Likewise, closed sets
  are \emph{down-sets} with respect to $\leq$.  If a space is
  Alexandroff, the converse holds as well: a set is open if and
    only if it is an up-set with respect to $\leq$, and it is closed
  if and only if it is a down-set.

  Now, let $\chi_{_{N}}$ be the conjunction of the following formulas (where
  we use $x<y$ as shorthand for $x\leq y\land x\neq y$):
  \[\begin{array}{l}
    \text{$\leq$ is a linear order} \\
    \ \ \ \forall xy.(x\leq y\ \land\ y\leq x\ \to\ x=y) \\
    \ \ \ \forall xy. (x\leq y\ \lor\ y\leq x) \\
    \\[-3mm]
    \text{There is a least element} \\
    \ \ \ \exists x.\forall y.(x\leq y) \\
    \\[-3mm]
    \text{Each element has an immediate successor in the ordering} \\
    \ \ \ \forall x.\exists y.(x < y\ \land\ \forall z.(x < z\ \to\ y\leq z)) \\
    \\[-3mm]
    \text{The space is Alexandroff (the down-sets are the closed sets)} \\
    \ \ \ \forall x.\exists U_x.\forall V_x.\forall y.(y\varepsilon V\to y\varepsilon U) \\
    \\[-3mm]
    \text{Each down-set other than $X, \emptyset$ has
       a least and a greatest element} \\
  \ \ \ \forall U.(\exists x.(x\nin U)\wedge \exists
x.(x\varepsilon U) \to\\
\ \ \ \ \ \ \ \ \ \ \ \ \ \ \ \ \ \ \exists z_l,z_g.(z_l\nin U\wedge
z_g\nin U\wedge \forall y.([y<z_l\vee z_g< y]\to y\varepsilon U)))
  \end{array}\]
  It is not hard to see that, if we take the open sets to be the
  up-sets, then $(\mathbb{N},\leq)$ is a model for $\chi_{_{N}}$.  In other
  words, $\chi_{_{N}}$ is satisfiable.  Now, suppose $\chi_{_{N}}$ is true in some
  topological space $(X,\tau)$.  We claim that $(X,\leq)$ is
  isomorphic to the natural numbers with their usual ordering. To
  prove this, it suffices to show that, for any $w\in X$, the set
  $\{v\mid v\leq w\}$ is finite (this property, together with the fact
  that $\leq$ is a linear order and each element has an immediate
  successor, characterizes the natural numbers up to isomorphism).
  In other words, we need to demonstrate that no infinite ascending or
  descending chains exist below an arbitrary point of $X$.

  Suppose that for some $w\in X$ the set $(w]=\set{v\in X\mid v\leq
    w}$ contains an infinite ascending chain $A{=}\set{a_1, a_2,\dots
  }$ with $a_i{<}a_{i+1}$ for each $i\in\N$. Consider the down-set
  $(A]=\set{w\in X\mid \exists a_i\in A.(w\leq a_i)}$ generated by the
  set $A$. Since $(A]\subseteq (w]$, we know that $(A]\neq X$, and
  hence there is a greatest element $g\in (A]$. By the definition of $(A]$,
   we have $g\leq a_i$ for some $i\in \N$. By definition of $A$, we
  also have $a_i<a_{i+1}$ and so $g<a_{i+1}$, contradicting the
  maximality of $g$ in $(A]$. Hence no infinite ascending chains exist
  below $w$.

  Next, suppose that for some $w\in X$ the set $(w]=\set{v\in X\mid
    v\leq w}$ contains an infinite descending chain
  $D=\set{d_1,d_2,\dots},$ with $d_{i+1}<d_i$ for each $i\in \N$. Then
  the set $X\setminus [D)=\set{w\in X\mid \forall d_i\in D.(w<d_i)}$
  is a non-empty down-set (for non-emptyness, note that the least
  element of $X$ cannot belong to $D$, and hence belongs to
  $X\setminus [D)$).  But then, there must be a greatest element
  $g\in X\setminus [D)$. Let $g'$ be the immediate successor of $g$.
  Note that by maximality of $g$ we must have $d_i\leq g'$ for some
  $i\in\N$. By definition of $D$, $d_{i+1}<d_i$ and we obtain
  $d_{i+1}\leq g$, hence $g\in [D)$, a contradiction. Thus no infinite
  descending chains exist below $w$.

\item Failure of Compactness

Consider the following set of $\mathcal L^2$-sentences with
one unary predicate $P$:
$$\Gamma \equiv \set{\chi_{_N}, \exists x.P(x)}\cup \set{\varphi_n\mid n\in\N}$$
where $\varphi_k\equiv \forall x.(P(x)\to\exists y_1,\dots,
y_k.(y_1<y_2<\ldots<y_k<x)$
express that every point in $P$ has at least $k$ predecessors.  Every
finite subset of $\Gamma$ is satisfiable but $\Gamma$ itself is not.

In fact, it is possible to show failure of compactness even without
using any unary predicates.

\item Failure of upward and downward L\"owenheim-Skolem

  Since $\chi_{_N}$ characterizes $(\mathbb{N},\leq)$ up to
  isomorphism, clearly, it has only countable models. Thus,
  the upward L\"owenheim-Skolem theorem fails for $\mathcal{L}^2$.
  The downwards L\"owenheim-Skolem theorem fails as well: we can
  easily express in $\mathcal L^2_=$ that the specialisation order
  $\leq$ is a dense linear ordering without endpoints. Further, we can
  express (on Alexandroff spaces) that each non-empty up-set has an
  infimum:
$$\text{Inf}\ \ \ \forall U.\big(\exists x.(x\varepsilon U)\to \exists
y.\forall z.\left((y<z\to z\varepsilon U)\wedge (z\varepsilon U\to
y\leq z)\right)\big)$$
Combining these formulas, we can enforce a \emph{complete dense linear
  order without endpoints}. An example of an infinite model satisfying
this is $\mathbb R$ with its usual ordering.  Any \emph{countable}
model, on the other hand, would have to be isomorphic to $\mathbb
Q$, as a countable dense linear order without endpoints, which
contradicts the conjunct Inf (e.g., the up-set $\set{w\in\mathbb
Q\mid w^2>2}$ has no infimum).
\item Failure of Interpolation

  Let $P, Q, R$ be distinct unary predicates. Let $\phi_{even}(P)$ be
  the $\mathcal{L}^2$-sentences expressing that, on the natural
  numbers, $P$ is true exactly of the even numbers, and
  $\phi_{even}(Q)$ likewise (it is not hard to see that there are such
  formulas). Then the following implication is valid:
  $$\chi_{_{N}}\land\phi_{even}(P)\land\exists x.(Px\land Rx)
     ~~~\to~~~ (\phi_{even}(Q)\to \exists x.(Qx\land Rx))$$
     Any interpolant for this implication has to express that $R$ is
     true of some even number, without the help of additional
     predicates. Using an Ehrenfeucht-Fra\"isse-style argument, one
     can show that this is impossible (note that we are essentially in
     first-order logic: quantification over open sets does provide
     any help, as the only open sets are the up-sets).

\item $\Sigma^1_1$-hard satisfiability problem.

  Using $\chi_{_{N}}$, we can reduce the problem of deciding
  whether an existential second order (ESO) formulas is true on
  $(\mathbb{N},\leq)$ ---a well known $\Sigma^1_1$-complete problem--- to
  the satisfiability problem of $\mathcal{L}^2$.  For simplity we will
  discuss here only the case for ESO sentences of the form $\exists
  R.\phi(R,{\leq})$, where $R$ is a single binary relation. The argument
  generalizes to more relations, and relations of other arities.

  Let an ESO sentence $\exists R.\phi$ be given. Let $N, P_1$ and
  $P_2$ be distinct unary predicates. Intuitively, the elements of the
  model satisfying $N$ will stand for natural numbers, while the other
  elements only play a technical role for coding up the binary
  relation $R$.  Let $x<^+y$ be short for $x<y\: \land\: Ny\: \land\:
  \forall z.(x<z\: \land\: Nz\: \to\: y\leq z)$, expressing that $y$ is the
  least $N$-element greater than $x$.  By induction, we define an
  $\mathcal{L}^2$-formula $\phi^*$ as follows:
  \[\begin{array}{lll}
     (x=y)^* &=& Nx \land Ny \land x=y \\
     (x\leq y)^* &=& Nx \land Ny \land x\leq y \\
     (Rxy)^* &=& \exists x'y'z.(z{<}x'{<^+}x\ \wedge\ z{<}y'{<^+}y\ \land\ P_1x'\ \land\ P_2 y')  \\
     (\phi\land\psi)^* &=& \phi^*\land\psi^* \\
     (\neg\phi)^* &=& \neg \phi^* \\
     (\exists x.\phi)^* &=& \exists x.(Nx\land\phi^*)
  \end{array}\]
  We claim that $(\mathbb{N},\leq) \models \exists R.\phi$ iff
  $\phi^*\land \chi_{_{N}}^\circ$ is satisfiable, where
  $\chi_{_{N}}^\circ$ is the relativisation of $\chi_{_{N}}$
  to $N$ (i.e., the formula obtained from $\chi_{_{N}}$ by
  relativising all quantifiers by $N$, thus expressing that the
  subspace defined by $N$ with its specialisation order is isomorphic
  to $(\mathbb{N},\leq)$).

  The difficult direction is \emph{left-to-right}. We give a rough
  sketch. Suppose that $(\mathbb{N},\leq) \models \exists R.\phi$. Let
  $\mathcal{R}\subseteq \mathbb{N}\times\mathbb{N}$ be a witnessing
  binary relation. Now, we define our model for $\phi^*\land
  \chi_{_{N}}^\circ$ as follows: the subspace defined by $N$ is
  simply the Alexandroff topology generated by
  $(\mathbb{N},\leq)$. For each pair $(m,n)\in \mathcal{R}$, we create
  three distinct $\neg N$-elements, $(m,n)_0$, $(m,n)_1$ and $(m,n)_2$.  Then
  we make sure that $m$ is the least $N$-successor of $(m,n)_1$ and $P_1$ holds
  at $(m,n)_1$,
  $n$ is the least $N$-successor of $(m,n)_2$ and $P_2$ holds
  at $(m,n)_2$, $(m,n)_0<(m,n)_1$ and $(m,n)_0<(m,n)_2$. In this way, we ensure that, for any pair
  of natural numbers $m,n$, $(m,n)\in\mathcal{R}$ iff the
  $\mathcal{L}^2$-formula $(Rxy)^*$ is true of $(m,n)$ in the
  constructed model. Once this observation is made, the claim becomes
  easy to prove.
\end{itemize}
\vspace{-7mm}
\end{proof}

\end{document}